\documentclass[12pt]{article}
\usepackage{amsmath,amsfonts,amssymb,amsthm}
\input amssym.def
\begin{document}
\def \Z{\Bbb Z}
\def \C{\Bbb C}
\def \R{\Bbb R}
\def \Q{\Bbb Q}
\def \N{\Bbb N}

\def \A{{\mathcal{A}}}
\def \D{{\mathcal{D}}}
\def \E{{\mathcal{E}}}
\def \E{{\mathcal{E}}}
\def \H{\mathcal{H}}
\def \S{{\mathcal{S}}}
\def \wt{{\rm wt}}
\def \tr{{\rm tr}}
\def \span{{\rm span}}
\def \Res{{\rm Res}}
\def \Der{{\rm Der}}
\def \End{{\rm End}}
\def \Ind {{\rm Ind}}
\def \Irr {{\rm Irr}}
\def \Aut{{\rm Aut}}
\def \GL{{\rm GL}}
\def \Hom{{\rm Hom}}
\def \mod{{\rm mod}}
\def \ann{{\rm Ann}}
\def \ad{{\rm ad}}
\def \rank{{\rm rank}\;}
\def \<{\langle}
\def \>{\rangle}

\def \g{{\frak{g}}}
\def \h{{\hbar}}
\def \k{{\frak{k}}}
\def \sl{{\frak{sl}}}
\def \gl{{\frak{gl}}}

\def \be{\begin{equation}\label}
\def \ee{\end{equation}}
\def \bex{\begin{example}\label}
\def \eex{\end{example}}
\def \bl{\begin{lem}\label}
\def \el{\end{lem}}
\def \bt{\begin{thm}\label}
\def \et{\end{thm}}
\def \bp{\begin{prop}\label}
\def \ep{\end{prop}}
\def \br{\begin{rem}\label}
\def \er{\end{rem}}
\def \bc{\begin{coro}\label}
\def \ec{\end{coro}}
\def \bd{\begin{de}\label}
\def \ed{\end{de}}

\def\NO{\mbox{\,$\circ\atop\circ$}\,}

\newcommand{\m}{\bf m}
\newcommand{\n}{\bf n}
\newcommand{\nno}{\nonumber}
\newcommand{\nord}{\mbox{\scriptsize ${\circ\atop\circ}$}}
\newtheorem{thm}{Theorem}[section]
\newtheorem{prop}[thm]{Proposition}
\newtheorem{coro}[thm]{Corollary}
\newtheorem{conj}[thm]{Conjecture}
\newtheorem{example}[thm]{Example}
\newtheorem{lem}[thm]{Lemma}
\newtheorem{rem}[thm]{Remark}
\newtheorem{de}[thm]{Definition}
\newtheorem{hy}[thm]{Hypothesis}
\makeatletter \@addtoreset{equation}{section}
\def\theequation{\thesection.\arabic{equation}}
\makeatother \makeatletter

\begin{center}
{\Large \bf Trigonometric Lie algebras, affine Lie algebras,  and vertex algebras}
\end{center}

\begin{center}
{Haisheng Li$^{a}$\footnote{Partially supported by
 China NSF grants (Nos.11471268, 11571391,11671247)},
Shaobin Tan$^{b}$\footnote{Partially supported by China NSF grants (Nos.11471268, 11531004)},
and Qing Wang$^{b}$\footnote{Partially supported by
 China NSF grants (Nos.11531004, 11622107), Natural Science Foundation of Fujian Province
(No. 2016J06002)}\\
$\mbox{}^{a}$Department of Mathematical Sciences\\
Rutgers University, Camden, NJ 08102, USA\\
$\mbox{}^{b}$School of Mathematical Sciences, Xiamen University,
Xiamen 361005, China}
\end{center}

\begin{abstract}
In this paper, we explore natural connections among trigonometric Lie
algebras, (general) affine Lie algebras, and vertex algebras. Among the main
results, we obtain a realization of trigonometric Lie algebras as what were called the covariant
algebras of the affine Lie algebra $\widehat{\A}$ of  Lie algebra $\A=\gl_{\infty}\oplus \gl_{\infty}$
with respect to certain automorphism groups. We then
prove that restricted modules of level $\ell$ for trigonometric
Lie algebras naturally correspond to equivariant quasi modules for the affine vertex
algebras $V_{\widehat{\A}}(\ell,0)$ (or $V_{\widehat{\A}}(2\ell,0)$).
Furthermore, we determine irreducible modules and equivariant quasi modules for simple vertex algebra
$L_{\widehat{\A}}(\ell,0)$ with $\ell$ a positive integer. In particular, we
 prove that every quasi-finite unitary highest weight (irreducible) module
 of level $\ell$ for type $A$ trigonometric Lie algebra gives rise to
 an irreducible equivariant quasi $L_{\widehat{\A}}(\ell,0)$-module.
\end{abstract}

\section{Introduction}

Trigonometric Lie algebras, of types $A$, $B$, $C$, and $D$,
are a family of infinite-dimensional Lie algebras (see  \cite{ffz}, \cite{flor}, \cite{gkl1}, \cite{gkl}).
For example, the rank $2$ trigonometric (sine) Lie algebra $\widehat{A}_{\hbar}$ (of type $A$)  with a real
parameter $\hbar$ is the Lie algebra with generators $A_{\alpha,m}$  for
$(\alpha,m)\in \Z^{2}$ and ${\bf c}$, a central
element, subject to relations
\begin{eqnarray*}
[A_{\alpha,m},A_{\beta,n}]=2i \sin
\hbar(m\beta-n\alpha)A_{\alpha+\beta,m+n}+m\delta_{\alpha+\beta,0}\delta_{m+n,0}{\bf c}
\end{eqnarray*}
for $\alpha,\beta,m,n\in \Z$.
In literature, the rank $2$ sine Lie algebra is also known as the quantum $2$-torus Lie
algebra. It was known (see \cite{hoppe}, \cite{flor}, \cite{gkl1,gkl}, \cite{kr1}) that the sine
Lie algebra has a canonical relation with $\overline{\gl_{\infty}}$ (the
formal completion with a central extension). Especially,
quasi-finite unitary highest weight (irreducible) modules were classified in  \cite{kr1} and
an explicit vertex operator realization was obtained in \cite{flor} and \cite{gkl}.
In this current paper, we study trigonometric Lie algebras with a different perspective, to explore
potential natural connections  with vertex algebras.

It has been known that (untwisted) affine Lie algebras have a natural association
with vertex algebras and modules (cf. \cite{fz}, \cite{li-local}), whereas twisted affine Lie algebras
are naturally associated  with twisted modules for vertex algebras (see \cite{flm}, \cite{li-twisted}).
Let $\g$ be any (possibly infinite-dimensional)  Lie algebra, equipped with a non-degenerate
symmetric invariant bilinear form $\<\cdot,\cdot\>$.
Then we have an affine Lie algebra $\widehat{\g}=\g\otimes \C[t,t^{-1}]\oplus \C {\bf k}$.
For any complex number $\ell$, we have a vertex algebra $V_{\widehat{\g}}(\ell,0)$
and its simple quotient vertex algebra $L_{\widehat{\g}}(\ell,0)$.
It is known that the category of $V_{\widehat{\g}}(\ell,0)$-modules
is canonically isomorphic to the category of restricted $\widehat{\g}$-modules of level $\ell$.
An important problem is to classify irreducible $L_{\widehat{\g}}(\ell,0)$-modules.
Every irreducible $L_{\widehat{\g}}(\ell,0)$-module is necessarily an irreducible
restricted $\widehat{\g}$-module of level $\ell$, but that is not all in general.

Suppose that $\g$ is a finite-dimensional simple Lie algebra with $\<\cdot,\cdot\>$
the killing form normalized such that the squared length of long roots equals $2$.
(In this case, the affine Lie algebra $\widehat{\g}$ is an affine Kac-Moody algebra.)
Let $\ell$ be a positive integer.
It was known (see \cite{fz}, \cite{dl}, \cite{li-local}) that
irreducible $L_{\widehat{\g}}(\ell,0)$-modules are exactly the
integrable highest weight $\widehat{\g}$-modules of level $\ell$, which are also the
unitary highest weight $\widehat{\g}$-modules of level $\ell$.

Note that the particular way to associate vertex algebras and their modules to affine Lie algebras works
for a wide variety of infinite-dimensional Lie algebras.
However, it no longer works  for trigonometric Lie algebras, including quantum $2$-torus Lie algebra.
With this as the main driving force,  a theory of what were called quasi modules
for vertex algebras was developed in  \cite{li-gamma}.
The notion of quasi module generalizes that of module for a vertex algebra in a certain natural way.
Indeed,  by using this theory a much wider variety of Lie algebras
can be associated with vertex algebras.
In order to study quasi modules more efficiently, the notion of vertex algebra was enhanced to
a notion of $\Gamma$-vertex algebra.
Later, a notion of equivariant quasi module
was introduced for a vertex $G$-algebra in \cite{li-tlie}.
It was proved in  \cite{li-twisted-quasi} that for a vertex operator algebra $V$
with a finite order automorphism $\sigma$,
the category of $\sigma$-twisted $V$-modules is canonically isomorphic to the category of
equivariant quasi $V$-modules with $G=\<\sigma\>$.
In view of this, the notion of quasi module can be viewed as a generalization of the notion of twisted module.

In this current paper, we first prove that trigonometric Lie algebras are isomorphic to what were called
the covariant algebras
of the affine Lie algebra $\widehat{\A}$ with respect to certain infinite order automorphism groups,
where $\A$ is a Lie algebra isomorphic to $\gl_{\infty}\oplus \gl_{\infty}$. Then we
give an  isomorphism between
the category of restricted modules for trigonometric Lie algebras of level $\ell$
and the category of equivariant quasi modules for $V_{\widehat{\A}}(\ell,0)$ for type $A$ and
for $V_{\widehat{\A}}(2\ell,0)$ for types $B$ and $D$.
Furthermore, we determine $L_{\widehat{\A}}(\ell,0)$-modules and
equivariant quasi $L_{\widehat{\A}}(\ell,0)$-modules with $\ell$ a positive integer.
In the establishment of the isomorphism between trigonometric Lie algebras and the covariant Lie algebras of
affine Lie algebras, the theory of equivariant quasi modules for vertex $\Gamma$-algebras plays a vital role.

Now, we mention some key ingredients in this paper.
Recall that associated to any Lie algebra $\g$ with a non-degenerate
symmetric invariant bilinear form $\<\cdot,\cdot\>$, we have an affine Lie algebra
$\widehat{\g}=\g\otimes \C[t,t^{-1}]\oplus \C {\bf k}$.
Let  $G$ be an automorphism group of $\g$, which preserves the bilinear form $\<\cdot,\cdot\>$,
satisfying the condition that for $a,b\in \g$,
$$[ga,b]=0\   \mbox{ and } \  \<ga,b\>=0$$
for all but finitely many $g\in G$. On the other hand, let $\chi$ be a linear character of $G$.
Let $G$ act on $\widehat{\g}$ by
$$g(a\otimes t^{m})=\chi(g)^{m} (ga\otimes t^{m})\ \ \ \mbox{ for }g\in
G,\ a\in \g,\ m\in \Z.$$
The {\em covariant algebra} of affine Lie algebra $\widehat{\g}$ with respect to $(G,\chi)$
is defined (see \cite{gkk}, \cite{li-tlie}) to be the Lie algebra
with $\widehat{\g}$ as the generating space and with defining relations
\begin{eqnarray*}
g(a\otimes t^{m})=a\otimes t^{m}\ \ \ \mbox{ for }g\in
G,\ a\in \g,\ m\in \Z
\end{eqnarray*}
and
\begin{eqnarray*}
[a\otimes t^{m}+\mu {\bf k},b\otimes t^{n}+\nu {\bf
k}]=\sum_{g\in G}\chi(g)^{m}\left([ga,b]\otimes
t^{m+n}+m\<ga,b\>\delta_{m+n,0}{\bf k}\right)
\end{eqnarray*}
for $a,b\in \g,\ m,n\in \Z,\ \mu,\nu\in \C$.
In case that $G$ is finite, the above covariant Lie algebra is isomorphic to the subalgebra of $G$-fixed points
in $\widehat{\g}$. In view of this,
covariant Lie algebras of affine Lie algebras can be viewed as generalizations of twisted affine Lie algebras.

We here need a particular Lie algebra $\g$.
Let $\gl_{\infty}$ denote the associative algebra of doubly infinite complex
matrices with finitely many nonzero entries. Fix
a symmetric invariant bilinear form $\<\cdot,\cdot\>$ by
$$\< E_{m,n},E_{p,q}\>={\rm tr}(E_{m,n}E_{p,q})=\delta_{m,q}\delta_{n,p}
\ \ \ \mbox{ for }m,n,p,q\in \Z.$$
Set $\A={\rm span} \{ E_{m,n}\ |\ m+n\in 2\Z\}\subset \gl_{\infty}$,
a subalgebra which is isomorphic to $\gl_{\infty}\oplus \gl_{\infty}$.
We then have an affine Lie algebra
$\widehat{\A}=\A\otimes \C[t,t^{-1}]\oplus \C {\bf k}.$
Let $\sigma$ be the automorphism of $\gl_{\infty}$ and $\A$ defined by
$\sigma(E_{m,n})=E_{m+1,n+1}$ for $m,n\in \Z$.
It is proved in this paper that $\widehat{A}_{\hbar}$ is isomorphic to
the covariant Lie algebra of the affine Lie algebra $\widehat{\A}$ with respect to automorphism group
$G=\<\sigma\>$ with linear character $\chi$ given by
$\chi(\sigma^k)=e^{ik\hbar}$ for $k\in \Z$.
Let $\tau$ be the linear endomorphism of $\gl_{\infty}$ defined by
$\tau(E_{m,n})=-E_{n,m}$ for $m,n\in \Z$. Then $\tau$ is an order-$2$ automorphism of $\gl_{\infty}$ and $\A$ viewed as Lie algebras. We prove that $\widehat{B}_{\hbar}$ and $\widehat{D}_{\hbar}$ are
isomorphic to the covariant Lie algebras of the affine Lie algebra $\widehat{\A}$ with respect to automorphism group
$\<\tau\>\times \<\sigma\>$ with suitably defined linear characters.
It is also proved that Lie algebras $\widehat{B}_{\hbar}$ and $\widehat{C}_{\hbar}$
are actually isomorphic.  Furthermore, we determine irreducible equivariant quasi modules
for simple vertex algebra
$L_{\widehat{\A}}(\ell,0)$.

For highest weight modules over an affine Kac-Moody algebra $\widehat{\g}$, it was known (see \cite{kac1}) that
integrability amounts to unitarity. On the other hand, integrability is equivalent to
 the condition that $e_{\theta}(x)^{\ell+1}=0$ (cf. \cite{lp}), where $e_{\theta}$ is a nonzero root vector of
the highest positive long root $\theta$. This equivalent condition plays a very important role
in the study of affine vertex operator algebras (cf. \cite{dl}, \cite{li-local}, \cite{mp1}, \cite{mp}).
Here, as a key result we establish an analogous condition
that is necessary and sufficient for a restricted $\widehat{A}_{\hbar}$-module of level $\ell$
to be an equivariant quasi $L_{\widehat{\A}}(\ell,0)$-module.
(This condition is similar to the quantum integrability condition studied in \cite{dm}, \cite{df}.)
Then, by using a result of Kac and Radul (see \cite{kr1}), we
 prove that every quasi-finite unitary highest weight (irreducible) $\widehat{A}_{\hbar}$-module
 of level $\ell$ gives rise to an
irreducible equivariant quasi $L_{\widehat{\A}}(\ell,0)$-module.
We conjecture that these are all irreducible equivariant quasi $L_{\widehat{\A}}(\ell,0)$-modules.

This paper is organized as follows: In Section 2, we give a realization of trigonometric Lie algebras
in terms of the covariant algebras of
affine Lie algebras of certain infinite-dimensional Lie algebras.
In Section 3, we establish an equivalence between
the category of restricted  $\widehat{A}_{\hbar}$-modules of a level $\ell$ and that of
equivariant quasi modules for vertex algebra $V_{\widehat{\A}}(\ell,0)$.
 In Section 4, we study modules for  simple vertex algebra $L_{\widehat{\A}}(\ell,0)$ with $\ell$ a positive integer.
In Section 5, we study equivariant quasi modules for $L_{\widehat{\A}}(\ell,0)$.

\section{Interplay between trigonometric Lie algebras and affine Lie algebras}
In this section, we first recall the trigonometric Lie algebras
$\widehat{A}_{\hbar}$, $\widehat{B}_{\hbar}$, $\widehat{C}_{\hbar}$,
$\widehat{D}_{\hbar}$, and then give an intrinsic connection of them with  the affine Lie algebra $\widehat{\A}$ of a
subalgebra $\A$ of the doubly infinite matrix Lie algebra $\gl_{\infty}$.

We begin by recalling the rank $2$ sine Lie algebra
(see \cite{ffz}, \cite{gkl1}).

\bd{dsinealgebra} {\em The rank $2$
sine Lie algebra $\widehat{A}_{\hbar}$ with a real
parameter $\hbar$ is generated by
$$A_{\alpha,m}\ \ \mbox{ for }
(\alpha,m)\in \Z^{2},$$ and by ${\bf c}$, a distinguished central
element, subject to relations
\begin{eqnarray}\label{edefining-relations}
[A_{\alpha,m},A_{\beta,n}]=2i \sin
\hbar(m\beta-n\alpha)A_{\alpha+\beta,m+n}+m\delta_{\alpha+\beta,0}\delta_{m+n,0}{\bf
c}
\end{eqnarray}
for $\alpha,\beta,m,n\in \Z$, where $i=\sqrt{-1}$.}
\ed

\br{ra00} {\em Note that  $A_{0,0}$ was excluded in the
original definition. The currently defined $\widehat{A}_{\hbar}$ is a (split)
direct sum of the originally defined $\widehat{A}_{\hbar}$ and $\C
A_{0,0}$. This slight modification is needed from a vertex-operator
point of view, just as with affine Heisenberg Lie algebras.  } \er

It can be readily seen that $\widehat{A}_{\hbar}$ is
a $\Z$-graded Lie algebra with
\begin{eqnarray}
\deg {\bf c}=0\   \mbox{ and }\deg A_{\alpha,m}=m\   \mbox{  for }\alpha,m\in \Z.
\end{eqnarray}
 We shall mainly use this $\Z$-grading in this paper.
(Note that $\widehat{A}_{\hbar}$ is also a $\Z^{2}$-graded
Lie algebra with $\deg A_{\alpha,m}=(\alpha,m)$ for $\alpha,m\in \Z$.)
The degree-zero subspace, which is linearly spanned by ${\bf
c}$ and $A_{\alpha,0}$ with $\alpha\in \Z$, is an abelian
subalgebra. Set
\begin{eqnarray}
H={\rm span}\{ A_{\alpha,0}\ |\ \alpha\in \Z\}.
\end{eqnarray}
The notion of $\Z$-graded $\widehat{A}_{\hbar}$-module is defined in the
obvious way. A $\Z$-graded $\widehat{A}_{\hbar}$-module is said to
be {\em quasifinite} (see \cite{kr1}) if every homogeneous subspace
is finite-dimensional.

An $\widehat{A}_{\hbar}$-module $W$ is said to be of {\em level}
$\ell\in \C$ if ${\bf c}$ acts on $W$ as scalar $\ell$. An
$\widehat{A}_{\hbar}$-module $W$ of level $\ell$ is called a {\em
highest weight module of highest weight $\lambda\in H^{*}$} if there
exists a vector $v\in W$ such that
\begin{eqnarray*}
&&A_{\alpha,n}v=0\ \ \ \ \ \mbox{ for all }\alpha,n\in \Z\
\mbox{with
}n>0,\\
&&hv=\lambda(h)v\ \ \ \mbox{ for }h\in H,\\
&&W=U(\widehat{A}_{\hbar})v.
\end{eqnarray*}
For $\ell\in \C,\ \lambda\in H^{*}$, denote by $M(\ell,\lambda)$ the
generalized Verma module. By definition, $M(\ell,\lambda)$ is the $\widehat{A}_{\hbar}$-module induced from
the $1$-dimensional module $\C$ for $\widehat{A}_{\hbar}^{\ge 0}$, on which $h$ acts as $\lambda(h)$ for $h\in H$,
${\bf c}$ acts as scalar $\ell$, and
$\widehat{A}_{\hbar}^{+}$ acts trivially.
Furthermore, denote by $L(\ell,\lambda)$ the irreducible
quotient module of $M(\ell,\lambda)$ by its (unique) maximal submodule.

Set
$$q=e^{i\hbar}\ \left(=e^{\hbar \sqrt{-1}}\right)\in \C^{\times}.$$
For $\alpha\in \Z$, form a generating function
\begin{eqnarray}
A_{\alpha}(z)=\sum_{n\in \Z}A_{\alpha,n}z^{-n-1}.
\end{eqnarray}
Then the defining relations (\ref{edefining-relations}) can be
written as
\begin{eqnarray}\label{exalpha-xbeta}
&&[A_{\alpha}(x),A_{\beta}(z)]
=q^{\alpha}A_{\alpha+\beta}(q^{\alpha}z)x^{-1}\delta\left(\frac{q^{\alpha+\beta}z}{x}\right)
\nonumber\\
&&\hspace{1cm}-
q^{-\alpha}A_{\alpha+\beta}(q^{-\alpha}z)x^{-1}\delta\left(\frac{q^{-(\alpha+\beta)}z}{x}\right)
+\delta_{\alpha+\beta,0}\frac{\partial}{\partial
z}x^{-1}\delta\left(\frac{z}{x}\right){\bf c}
\end{eqnarray}
for $\alpha,\beta\in \Z$.

Next, we relate $\widehat{A}_{\hbar}$ to an affine Lie algebra. Let $\g$
be a (possibly infinite-dimensional) Lie algebra equipped with a
non-degenerate symmetric invariant bilinear form $\<\cdot,\cdot\>$.
Then one has an {\em (untwisted) affine Lie algebra}
$$\hat{\g}=\g\otimes \C[t,t^{-1}]\oplus \C {\bf k},$$
where ${\bf k}$ is a nonzero central element and
\begin{eqnarray}\label{eaffine-defining-relation-1}
[a\otimes t^{m},b\otimes t^{n}]=[a,b]\otimes
t^{m+n}+m\delta_{m+n,0}\<a,b\>{\bf k}
\end{eqnarray}
for $a,b\in \g,\ m,n\in \Z$.

The following was proved in \cite{li-tlie} (cf. \cite{gkk}):

\bp{ptwisted-Lie} Let $\g$ be a Lie algebra equipped with a
symmetric invariant bilinear form $\<\cdot,\cdot\>$, let $\Gamma$ be
a subgroup of ${\rm Aut}(\g,\<\cdot,\cdot\>)$, and let $\phi:
\Gamma\rightarrow \C^{\times}$ be a linear character (group
homomorphism). Assume that for $a,b\in \g$,
\begin{eqnarray}
[ga,b]=0\ \ \mbox{ and }\ \ \<ga,b\>=0\ \ \mbox{ for all but
finitely many }g\in \Gamma.
\end{eqnarray}
Define a bilinear operation $[\cdot,\cdot]_{\Gamma}$
on the vector space $\g\otimes \C[t,t^{-1}]\oplus \C {\bf k}$ by
\begin{eqnarray}
[a\otimes t^{m}+\mu {\bf k},b\otimes t^{n}+\nu {\bf
k}]_{\Gamma}=\sum_{g\in \Gamma}\phi(g)^{m}\left([ga,b]\otimes
t^{m+n}+m\<ga,b\>\delta_{m+n,0}{\bf k}\right)
\end{eqnarray}
for $a,b\in \g,\ m,n\in \Z,\ \mu,\nu\in \C$. Then the subspace
linearly spanned by elements
\begin{eqnarray}
\phi(g)^{m}(ga\otimes t^{m})-(a\otimes t^{m})\ \ \ \mbox{ for }g\in
\Gamma,\ a\in \g,\ m\in \Z
\end{eqnarray}
is a two-sided ideal of the nonassociative algebra $(\g\otimes
\C[t,t^{-1}]+ \C {\bf k}, [\cdot,\cdot]_{\Gamma})$, and the quotient algebra, which is denoted by
$\hat{\g}[\Gamma]$, is a Lie algebra. \ep

The Lie algebra $\hat{\g}[\Gamma]$ obtained in Proposition
\ref{ptwisted-Lie} is called the {\em $(\Gamma,\phi)$-covariant
algebra of the affine Lie algebra $\hat{\g}$}, or simply the {\em $\Gamma$-covariant algebra of
$\hat{\g}$} whenever the context is clear. For $a\in \g,\ m\in \Z$,
we denote by $\overline{a\otimes t^{m}}$ the image of $a\otimes
t^{m}$ under the natural map from $\hat{\g}$ onto
$\hat{\g}[\Gamma]$. We shall still use ${\bf k}$ for its image.

Let $\gl_{\infty}$ be the associative algebra of all doubly infinite
complex matrices with only finitely many nonzero entries. Naturally,
$\gl_{\infty}$ is a Lie algebra under the commutator bracket. For
$m,n\in \Z$, let $E_{m,n}$ denote the unit matrix whose only nonzero
entry is the $(m,n)$-entry which is $1$. Then $E_{m,n}$ ($m,n\in
\Z$) form a basis of $\gl_{\infty}$, where
\begin{eqnarray}
E_{m,n}\cdot E_{p,q}=\delta_{n,p}E_{m,q},\ \ \ \
[E_{m,n},E_{p,q}]=\delta_{n,p}E_{m,q}-\delta_{q,m}E_{p,n}
\end{eqnarray}
for $m,n,p,q\in \Z$. Equip $\gl_{\infty}$ with the bilinear form
$\<\cdot,\cdot\>$ defined by
\begin{eqnarray}
\< E_{m,n},E_{r,s}\>={\rm
tr}(E_{m,n}E_{r,s})=\delta_{m,s}\delta_{n,r}
\end{eqnarray}
for $m,n,r,s\in \Z$. This bilinear form is non-degenerate,
symmetric, and associative.

\bd{dsigma-r}{\em For $r\in \Z$, define a linear operator
 $\sigma_{r}$ on $\gl_{\infty}$ by
\begin{eqnarray}
\sigma_{r}(E_{m,n})=E_{m+r,n+r}\ \ \ \mbox{ for }m,n\in \Z.
\end{eqnarray}}
\ed

The following is straightforward to prove (cf. \cite{gkk}):

\bl{lZ-action} The map $\Z\ni n\mapsto \sigma_{n}$ defines a group action
of $\Z$ on $\gl_{\infty}$ by automorphisms which preserve the
bilinear form $\<\cdot,\cdot\>$. Furthermore, for any $a,b\in
\gl_{\infty}$, we have
\begin{eqnarray}
[\sigma_{r}(a),b]=0\ \ \mbox{ and }\ \ \<\sigma_{r}(a),b\>=0
\end{eqnarray}
for all but finitely many $r\in \Z$. \el

\br{rnew-algebra}  {\em As a key ingredient, we here introduce an
associative algebra. Set
\begin{eqnarray}
{\mathcal{A}}={\rm span} \{E_{m,n}\ | \ m,n\in \Z\  \mbox{ with
}m+n\in 2\Z\},
\end{eqnarray}
which is an associative subalgebra of $\gl_{\infty}$. For $\alpha,m\in \Z$, set
\begin{eqnarray}\label{edefG}
G_{\alpha,m}=E_{\alpha+m,m-\alpha}\in \A.
\end{eqnarray}
Then $G_{\alpha,m}$ $(\alpha,m\in \Z)$ form a basis of ${\mathcal{A}}$ and we have
\begin{eqnarray}\label{eGalpham}
&&G_{\alpha,m}\cdot
G_{\beta,n}=\delta_{m-\alpha,n+\beta}G_{\alpha+\beta,\alpha+n},\nonumber\\
&&\<G_{\alpha,m},G_{\beta,n}\>=\delta_{\alpha+\beta,0}\delta_{m,n}.
\end{eqnarray}
Furthermore, it can be readily seen that $\A$ is stable under the
action of $\Z$, where
\begin{eqnarray}
\sigma_{r}(G_{\alpha,m})=G_{\alpha,m+r} \ \ \ \mbox{ for
}r,\alpha,m\in \Z.
\end{eqnarray}}
\er

View $\A$ as a Lie algebra and equip $\A$ with the symmetric
invariant bilinear form $\<\cdot,\cdot\>$ defined above, which is
still non-degenerate. Then we have an affine Lie algebra
\begin{eqnarray}
\widehat{\A}=\A\otimes \C[t,t^{-1}]\oplus \C {\bf k},
\end{eqnarray}
on which $\Z$ acts as an automorphism group.

Define a linear character $\chi_{q}:\Z\rightarrow \C^{\times}$ by
\begin{eqnarray}
\chi_{q}(n)=q^{n}\ \ \ \mbox{ for }n\in \Z,
\end{eqnarray}
recalling that $q=e^{i\hbar}$ where $\hbar$ is the same one as for
$\widehat{A}_{\hbar}$. Then we have (cf. \cite{gkk}):

\bp{ptypeA} The sine Lie algebra $\widehat{A}_{\hbar}$ is isomorphic to
the $(\Z,\chi_{q})$-covariant algebra $\widehat{\A}[\Z]$ of the
affine Lie algebra $\widehat{\A}$ with ${\bf c}={\bf k}$ and with
\begin{eqnarray}
A_{\alpha,m}=\overline{G_{\alpha,0}\otimes t^{m}} \ \ \ \mbox{ for
}\alpha,m\in \Z.
\end{eqnarray}
 \ep

\begin{proof} {}From the definition, $\widehat{A}_{\hbar}$ has a
basis $\{ A_{\alpha,m}\ |\ \alpha,m\in \Z\}\cup \{ \bf c\}$ with the
given bracket relations. On the other hand, {}from Proposition
\ref{ptwisted-Lie},
 $\widehat{\A}[\Z]$ as a vector space is the quotient space of
$\widehat{\A}$, modulo the subspace linearly spanned by
$$\sigma_{n}(a)\otimes t^{m}-q^{-mn}(a\otimes t^{m})$$
for $a\in \A,\; m,n\in \Z$, and its bracket relation is given by
\begin{eqnarray}\label{ebracket-relation}
[\overline{a\otimes t^{m}},\overline{b\otimes t^{n}}]=\sum_{r\in
\Z}q^{mr}\left(\overline{[\sigma_{r}(a),b]\otimes
t^{m+n}}+m\<\sigma_{r}(a),b\>\delta_{m+n,0} {\bf k}\right),
\end{eqnarray}
where  $\overline{u\otimes t^{p}}$ denotes the image of $u\otimes
t^{p}$ in $\widehat{\A}[\Z]$ for $u\in \A,\; p\in \Z$, and ${\bf k}$
is identified with its image. We see that ${\bf k}$ and
$\overline{G_{\alpha,0}\otimes t^{m}}$ ($\alpha,m\in \Z$) form a
basis of $\widehat{\A}[\Z]$.
Let $\alpha,\beta,m,n\in \Z$.  We have
\begin{eqnarray}\label{esub-calculation}
&&[\sigma_{r}(G_{\alpha,0}),G_{\beta,0}]=[G_{\alpha,r},G_{\beta,0}]
=\delta_{r,\alpha+\beta}G_{\alpha+\beta,\alpha}-\delta_{r,-\alpha-\beta}G_{\alpha+\beta,-\alpha},
\nonumber\\
&&
\<\sigma_{r}(G_{\alpha,0}),G_{\beta,0}\>=\<G_{\alpha,r},G_{\beta,0}\>=\delta_{\alpha+\beta,0}\delta_{r,0}
\end{eqnarray}
for $r\in \Z$. Applying this to (\ref{ebracket-relation}) we get
\begin{eqnarray}
&&[\overline{G_{\alpha,0}\otimes t^{m}},\overline{G_{\beta,0}\otimes
t^{n}}]\nonumber\\
&=&q^{m(\alpha+\beta)}\overline{G_{\alpha+\beta,\alpha}\otimes
t^{m+n}}-q^{-m(\alpha+\beta)}\overline{G_{\alpha+\beta,-\alpha}\otimes
t^{m+n}}+m\delta_{\alpha+\beta,0}\delta_{m+n,0}{\bf k}\nonumber\\
&=&q^{m\beta-n\alpha}\overline{G_{\alpha+\beta,0}\otimes
t^{m+n}}-q^{n\alpha-m\beta}\overline{G_{\alpha+\beta,0}\otimes
t^{m+n}}+m\delta_{\alpha+\beta,0}\delta_{m+n,0} {\bf k}\nonumber\\
&=&\left(q^{m\beta-n\alpha}-q^{n\alpha-m\beta}\right)\overline{G_{\alpha+\beta,0}\otimes
t^{m+n}}+m\delta_{\alpha+\beta,0}\delta_{m+n,0} {\bf k},
\end{eqnarray}
noticing that
\begin{eqnarray*}
&&\overline{G_{\alpha+\beta,\alpha}\otimes t^{m+n}}
=\overline{\sigma_{\alpha}G_{\alpha+\beta,0}\otimes t^{m+n}}
=q^{-(m+n)\alpha}\overline{G_{\alpha+\beta,0}\otimes t^{m+n}},\\
&&\overline{G_{\alpha+\beta,-\alpha}\otimes t^{m+n}}
=\overline{\sigma_{-\alpha}G_{\alpha+\beta,0}\otimes t^{m+n}}
=q^{(m+n)\alpha}\overline{G_{\alpha+\beta,0}\otimes t^{m+n}}.
\end{eqnarray*}
Now it follows that $\widehat{A}_{\hbar}$ is isomorphic to
$\widehat{\A}[\Z]$ with $A_{\alpha,m}$ corresponding to
$\overline{G_{\alpha,0}\otimes t^{m}}$ for $\alpha,m\in \Z$ and with
${\bf c}$ corresponding to ${\bf k}$.
\end{proof}

Next, we discuss Lie algebra $\widehat{B}_{\hbar}$.
Define a second order automorphism $\tau_{B}$ of $\widehat{A}_{\hbar}$ by
\begin{eqnarray}
\tau_{B}(A_{\alpha,m})=-(-1)^{m}A_{-\alpha,m}\   \    \mbox{ for }\alpha,m\in \Z.
\end{eqnarray}
Lie algebra $\widehat{B}_{\hbar}$ is defined to be the subalgebra
$(\widehat{A}_{\hbar})^{\tau_{B}}$ of $\tau_{B}$-fixed points
in $\widehat{A}_{\hbar}$, which is linearly spanned by elements
\begin{eqnarray}
B_{\alpha,m}:=A_{\alpha,m}-(-1)^{m}A_{-\alpha,m}
\end{eqnarray}
for $\alpha,m\in \Z$.
The commutation relations are
\begin{eqnarray}\label{eB-relations}
&&[B_{\alpha,m},B_{\beta,n}]=2i\sin\hbar(m\beta-n\alpha)B_{\alpha+\beta,m+n}+(-1)^{n}2i\sin
\hbar(m\beta+n\alpha)B_{\alpha-\beta,m+n}\nonumber\\
&&\hspace{2cm}
+2m\left(\delta_{\alpha+\beta,0}-(-1)^{m}\delta_{\alpha-\beta,0}\right)\delta_{m+n,0}{\bf
c}.
\end{eqnarray}

For $\alpha\in \Z$, form a generating function
\begin{eqnarray}
B_{\alpha}(z)=\sum_{m\in \Z}B_{\alpha,m}z^{-m-1}.
\end{eqnarray}
That is,
\begin{eqnarray}
B_{\alpha}(z)=A_{\alpha}(z)+A_{-\alpha}(-z).
\end{eqnarray}
We have
\begin{eqnarray}\label{B-relation1}
B_{-\alpha}(z)=B_{\alpha}(-z)
\end{eqnarray}
for $\alpha\in \Z$, while (\ref{eB-relations})  in terms of generating functions becomes
\begin{eqnarray}\label{B-relation2}
&&[B_{\alpha}(x),B_{\beta}(z)]
=q^{\alpha}B_{\alpha+\beta}(q^{\alpha}z)x^{-1}\delta\left(\frac{q^{\alpha+\beta}z}{x}\right)
-q^{-\alpha}B_{\alpha+\beta}(q^{-\alpha}z)x^{-1}\delta\left(\frac{q^{-(\alpha+\beta)}z}{x}\right)\nonumber\\
&&\ \ \ \
+q^{\alpha}B_{\alpha-\beta}(-q^{\alpha}z)x^{-1}\delta\left(\frac{-q^{\alpha-\beta}z}{x}\right)
-q^{-\alpha}B_{\alpha-\beta}(-q^{-\alpha}z)x^{-1}\delta\left(\frac{-q^{-(\alpha-\beta)}z}{x}\right)\nonumber\\
&&\ \ \ \ +2\delta_{\alpha+\beta,\bar{0}}\frac{\partial}{\partial
z}x^{-1}\delta\left(\frac{z}{x}\right){\bf c}
-2\delta_{\alpha-\beta,0}\frac{\partial}{\partial
z}x^{-1}\delta\left(\frac{z}{-x}\right){\bf c}.
\end{eqnarray}

\bd{dtau} {\em Let $\tau$ be the order-2 automorphism of Lie algebra $\gl_{\infty}$ defined by
\begin{eqnarray}
\tau (E_{m,n})=-E_{n,m}\ \ \ \mbox{ for }m,n\in \Z.
\end{eqnarray}}
\ed

We summarize some straightforward facts as follows:

\bl{labout-tau} The automorphism $\tau$ commutes with the action of
$\Z$ on $\gl_{\infty}$ and also preserves the bilinear form
$\<\cdot,\cdot\>$.  The group $\Z_{2}\times \Z$ acts on
$\gl_{\infty}$, extending the action of $\Z$ defined in Lemma
\ref{lZ-action}, such that
$$\tau(E_{m,n})=-E_{n,m},\ \ \ \sigma_{r}(E_{m,n})=E_{m+r,n+r}
\ \ \ \mbox{ for }m,n,r\in \Z.$$ Furthermore, the subalgebra $\A$
 is stable under the action of
 $\Z_{2}\times \Z$ and we have
\begin{eqnarray}
\tau(G_{\alpha,m})=-G_{-\alpha,m}, \ \ \ \
\sigma_{r}(G_{\alpha,m})=G_{\alpha,m+r}\ \ \ \mbox{ for
}\alpha,m,r\in \Z.
\end{eqnarray}
 \el

Define a linear character $\chi_{q}^{B}: \Z_{2}\times \Z\rightarrow
\C^{\times}$ by
\begin{eqnarray}
\chi_{q}^{B}(\tau)=-1 \ \mbox{ and } \
\chi_{q}^{B}(\sigma_{r})=q^{r} \ \mbox{ for }r\in \Z,
\end{eqnarray}
where it is understood that $\Z_{2}=\<\tau\>$ and $\Z=\{ \sigma_{r}\
|\ r\in \Z\}$. Then we have:

\bp{ptypeB} Lie algebra $\widehat{B}_{\hbar}$ is isomorphic to the
$(\Z_{2}\times \Z,\chi_{q}^{B})$-covariant algebra
$\widehat{\A}[\Z_{2}\times\Z]$ of the affine Lie algebra
$\widehat{\A}$ with ${\bf c}=\frac{1}{2}{\bf k}$ and with
\begin{eqnarray}
B_{\alpha,m}=\overline{G_{\alpha,0}\otimes t^{m}} \ \ \ \mbox{ for
}\alpha,m\in \Z.
\end{eqnarray}
\ep

\begin{proof} {}From the definition, $\widehat{B}_{\hbar}$ is the Lie algebra with generators
$B_{\alpha,m}$ for $\alpha,m\in \Z$, subject to relations
\begin{eqnarray}
 B_{-\alpha,m}=-(-1)^{m}B_{\alpha,m},
\end{eqnarray}
\begin{eqnarray*}
&&[B_{\alpha,m},B_{\beta,n}]=\left(q^{m\beta-n\alpha}-q^{n\alpha-m\beta}\right)B_{\alpha+\beta,m+n}
\nonumber\\
&&\ \ \ \
+(-1)^{n}\left(q^{m\beta+n\alpha}-q^{-m\beta-n\alpha}\right)B_{\alpha-\beta,m+n}
+2m(\delta_{\alpha+\beta,0}-(-1)^{m}\delta_{\alpha-\beta,0})\delta_{m+n,0}{\bf
c}\hspace{0.5cm}
\end{eqnarray*}
for $\alpha,\beta,m,n\in \Z$.
On the other hand, consider the covariant algebra
$\widehat{\A}[\Z_{2}\times \Z]$ with respect to the linear character
$\chi_{q}^{B}$. Let $\alpha,\beta,m,n\in \Z$. From the definition,
we have
\begin{eqnarray*}
&&[\overline{G_{\alpha,0}\otimes t^{m}},\overline{G_{\beta,0}\otimes
t^{n}}]\nonumber\\
&=&\sum_{r\in
\Z}q^{mr}\left(\overline{[\sigma_{r}(G_{\alpha,0}),G_{\beta,0}]\otimes
t^{m+n}}+m\<\sigma_{r}(G_{\alpha,0}),G_{\beta,0}\>\delta_{m+n,0}
{\bf k}\right)\nonumber\\
&&+\sum_{r\in \Z}(-q^{r})^{m}\left(\overline{[\tau
\sigma_{r}(G_{\alpha,0}),G_{\beta,0}]\otimes
t^{m+n}}+m\<\tau\sigma_{r}(G_{\alpha,0}),G_{\beta,0}\>\delta_{m+n,0}
{\bf k}\right)\\
&=&\sum_{r\in
\Z}q^{mr}\left(\overline{[G_{\alpha,r},G_{\beta,0}]\otimes
t^{m+n}}+m\<G_{\alpha,r},G_{\beta,0}\>\delta_{m+n,0}
{\bf k}\right)\nonumber\\
&&-\sum_{r\in
\Z}(-q^{r})^{m}\left(\overline{[G_{-\alpha,r},G_{\beta,0}]\otimes
t^{m+n}}+m\<G_{-\alpha,r},G_{\beta,0}\>\delta_{m+n,0}
{\bf k}\right)\\
&=&\sum_{r\in \Z}q^{mr}\overline{[G_{\alpha,r},G_{\beta,0}]\otimes
t^{m+n}}- \sum_{r\in
\Z}(-q^{r})^{m}\overline{[G_{-\alpha,r},G_{\beta,0}]\otimes
t^{m+n}}\\
&&+ \sum_{r\in \Z}m\left(q^{mr}\delta_{\alpha+\beta,0}\delta_{r,0}
-(-q^{r})^{m}\delta_{\beta-\alpha,0}\delta_{r,0}\right)\delta_{m+n,0}
{\bf k}\\
&=&\sum_{r\in \Z}q^{mr}\overline{[G_{\alpha,r},G_{\beta,0}]\otimes
t^{m+n}}- \sum_{r\in
\Z}(-q^{r})^{m}\overline{[G_{-\alpha,r},G_{\beta,0}]\otimes
t^{m+n}}\\
&&+ m\left(\delta_{\alpha+\beta,0}
-(-1)^{m}\delta_{\alpha-\beta,0}\right)\delta_{m+n,0}{\bf k}.
\end{eqnarray*}
Just as with $\widehat{A}_{\hbar}$ we have
\begin{eqnarray*}
\sum_{r\in \Z}q^{mr}\overline{[G_{\alpha,r},G_{\beta,0}]\otimes
t^{m+n}}
=\left(q^{m\beta-n\alpha}-q^{n\alpha-m\beta}\right)\overline{G_{\alpha+\beta,0}\otimes
t^{m+n}},
\end{eqnarray*}
and we also have
\begin{eqnarray*}
&&\sum_{r\in \Z}(-q^{r})^{m}\overline{[G_{-\alpha,r},G_{\beta,0}]\otimes t^{m+n}}\\
&=&\sum_{r\in
\Z}(-q^{r})^{m}\left(\delta_{r,\beta-\alpha}\overline{G_{\beta-\alpha,-\alpha}\otimes
t^{m+n}}-\delta_{r,\alpha-\beta}\overline{G_{\beta-\alpha,\alpha}\otimes
t^{m+n}}\right)\\
&=&\sum_{r\in
\Z}(-q^{r})^{m}\left(\delta_{r,\alpha-\beta}\overline{\tau\sigma_{\alpha}G_{\alpha-\beta,0}\otimes
t^{m+n}}-\delta_{r,\beta-\alpha}\overline{\tau\sigma_{-\alpha}G_{\alpha-\beta,0}\otimes
t^{m+n}}\right)\\
&=&\sum_{r\in
\Z}\delta_{r,\alpha-\beta}(-q^{r})^{m}(-q^{\alpha})^{-m-n}\overline{G_{\alpha-\beta,0}\otimes
t^{m+n}}\\
&&-\sum_{r\in
\Z}\delta_{r,\beta-\alpha}(-q^{r})^{m}(-q^{-\alpha})^{-m-n}\overline{G_{\alpha-\beta,0}\otimes
t^{m+n}}\\
&=&(-1)^{n}\left(q^{-m\beta-n\alpha}-q^{m\beta+n\alpha}\right)\overline{G_{\alpha-\beta,0}\otimes
t^{m+n}}.
\end{eqnarray*}
Therefore
\begin{eqnarray}
&&[\overline{G_{\alpha,0}\otimes t^{m}},\overline{G_{\beta,0}\otimes
t^{n}}]\nonumber\\
&=&\left(q^{m\beta-n\alpha}-q^{n\alpha-m\beta}\right)\overline{G_{\alpha+\beta,0}\otimes
t^{m+n}}+
(-1)^{n}\left(q^{m\beta+n\alpha}-q^{-m\beta-n\alpha}\right)\overline{G_{\alpha-\beta,0}\otimes
t^{m+n}}\nonumber\\
&&+m\left(\delta_{\alpha+\beta,0}
-(-1)^{m}\delta_{\alpha-\beta,0}\right)\delta_{m+n,0}{\bf k}.
\end{eqnarray}
On the other hand, for $\alpha,m\in \Z$ we have
\begin{eqnarray*}
\overline{G_{-\alpha,0}\otimes t^{m}}=-\overline{\tau
G_{\alpha,0}\otimes
t^{m}}=-\chi(\tau)^{-m}\overline{G_{\alpha,0}\otimes t^{m}}
=-(-1)^{m}\overline{G_{\alpha,0}\otimes t^{m}}.
\end{eqnarray*}
Therefore, $\widehat{B}_{\hbar}$ is isomorphic to
$\widehat{\A}[\Z_{2}\times \Z]$ with $B_{\alpha,m}$ corresponding to
$\overline{G_{\alpha,0}\otimes t^{m}}$ for $\alpha,m\in \Z$ and with
${\bf c}$ corresponding to $\frac{1}{2}{\bf k}$.
\end{proof}

The Lie algebra $\widehat{C}_{\hbar}$ is the subalgebra of $\tau_{C}$-fixed points in
$\widehat{A}_{\hbar}$, where $\tau_{C}$ is the second order automorphism  of $\widehat{A}_{\hbar}$,
defined by
\begin{eqnarray}
\tau_{C}(A_{\alpha,m})=-(-1)^{m}q^{2\alpha}A_{-\alpha,m}\ \ \ \mbox{
for }\alpha,m\in \Z.
\end{eqnarray}
As a subspace,  $\widehat{C}_{\hbar}$ is linearly spanned by
\begin{eqnarray}
C_{\alpha,m}:=A_{\alpha,m}-(-1)^{m}q^{2\alpha}A_{-\alpha,m}\ \ \
\mbox{ for }\alpha,m\in \Z,
\end{eqnarray}
where the commutation relations  are
\begin{eqnarray}
[C_{\alpha,m},C_{\beta,n}]&=&2i\sin\hbar(m\beta-n\alpha)C_{\alpha+\beta,m+n}\nonumber\\
&&+(-1)^{n}2iq^{2\beta}\sin\hbar
(m\beta+n\alpha)C_{\alpha-\beta,m+n}\nonumber\\
&&+2m\left(\delta_{\alpha+\beta,0}-(-1)^{m}q^{2\alpha}\delta_{\alpha-\beta,0}\right)\delta_{m+n,0}{\bf
c}.\ \hspace{1cm}
\end{eqnarray}
In terms of generating functions we have
\begin{eqnarray}
C_{\alpha}(x)=A_{\alpha}(x)+q^{2\alpha}A_{-\alpha}(-x).
\end{eqnarray}
Note that
$$C_{-\alpha}(x)=A_{-\alpha}(x)+q^{-2\alpha}A_{\alpha}(-x)=q^{-2\alpha}C_{\alpha}(-x).$$
We have:

\bp{ptypeC} Lie algebra $\widehat{C}_{\hbar}$ is isomorphic to $\widehat{B}_{\hbar}$,
which  is isomorphic to the
$(\Z_{2}\times \Z,\chi_{q}^{B})$-covariant algebra
$\widehat{\A}[\Z_{2}\times \Z]$.
\ep

\begin{proof} For $\alpha\in \Z$, set
\begin{eqnarray}
C'_{\alpha}(x)=q^{-\alpha}C_{\alpha}(x)=q^{-\alpha}A_{\alpha}(x)+q^{\alpha}A_{-\alpha}(-x).
\end{eqnarray}
Then we have
$$C'_{-\alpha}(x)=C'_{\alpha}(-x)$$
and
\begin{eqnarray}
[C_{\alpha,m}',C_{\beta,n}']&=&2i\sin \hbar(m\beta-n\alpha)C_{\alpha+\beta,m+n}'\nonumber\\
&&+(-1)^{n}2i\sin \hbar
(m\beta+n\alpha)C'_{\alpha-\beta,m+n}\nonumber\\
&&+2m\left(\delta_{\alpha+\beta,0}-(-1)^{m}\delta_{\alpha-\beta,0}\right)\delta_{m+n,0}{\bf
c}. \ \hspace{1cm}
\end{eqnarray}
Comparing these with (\ref{B-relation1}) and (\ref{B-relation2}),
we conclude that $\widehat{C}_{\hbar}$ is isomorphic to $\widehat{B}_{\hbar}$
with $C'_{\alpha}(x)$ corresponding to $B_{\alpha}(x)$ for $\alpha\in \Z$ and with ${\bf c}$
corresponding to ${\bf c}$.
\end{proof}

Last, we discuss the Lie algebra $\widehat{D}_{\hbar}$.
Let $\tau_{D}$ be the linear endomorphism of $\widehat{A}_{\hbar}$ defined by
\begin{eqnarray}
\tau_{D}({\bf c})={\bf c}, \   \  \  \
 \tau_{D}(A_{\alpha,m})=-q^{2\alpha}A_{-\alpha,m}\ \ \ \mbox{ for }\alpha,m\in \Z.
\end{eqnarray}
It is straightforward to show that $\tau_{D}$ is an
order-$2$ automorphism  of $\widehat{A}_{\hbar}$.
Lie algebra $\widehat{D}_{\hbar}$ is defined to be the $\tau_{D}$-fixed points subalgebra of $A_{\hbar}$.

For $\alpha,m\in \Z$, set
\begin{eqnarray}
D_{\alpha,m}=A_{\alpha,m}-q^{2\alpha}A_{-\alpha,m}.
\end{eqnarray}
Then $\widehat{D}_{\hbar}$ is linearly spanned by $D_{\alpha,m}$ for $\alpha,m\in \Z$.
We have
$$D_{-\alpha,m}=-q^{-2\alpha}D_{\alpha,m}\   \    \    \mbox{ for }\alpha,m\in \Z$$
and
\begin{eqnarray}
&&[D_{\alpha,m},D_{\beta,n}]=2i\sin\hbar(m\beta-n\alpha)D_{\alpha+\beta,m+n}\nonumber\\
&&\ \ \ \ +2iq^{2\beta}\sin
\hbar(m\beta+n\alpha)D_{\alpha-\beta,m+n}
+2m(\delta_{\alpha+\beta,0}-q^{2\alpha}\delta_{\alpha-\beta,0})\delta_{m+n,0}\cdot
{\bf c}.\ \hspace{1cm}
\end{eqnarray}

For $\alpha,m\in \Z$, set
\begin{eqnarray}
D_{\alpha}(x)=\sum_{m\in \Z}D_{\alpha,m}x^{-m-1}.
\end{eqnarray}
In terms of generating functions, we have
\begin{eqnarray}
&&[D_{\alpha}(x),D_{\beta}(z)]
=q^{\alpha}D_{\alpha+\beta}(q^{\alpha}z)x^{-1}\delta\left(\frac{q^{\alpha+\beta}z}{x}\right)
-q^{-\alpha}D_{\alpha+\beta}(q^{-\alpha}z)x^{-1}\delta\left(\frac{q^{-(\alpha+\beta)}z}{x}\right)\nonumber\\
&&\ \ \ \ \ \
+q^{2\beta-\alpha}D_{\alpha-\beta}(q^{-\alpha}z)x^{-1}\delta\left(\frac{q^{-(\alpha-\beta)}z}{x}\right)
-q^{\alpha+2\beta}D_{\alpha-\beta}(q^{\alpha}z)x^{-1}\delta\left(\frac{q^{\alpha-\beta}z}{x}\right)
\nonumber\\
&&\ \ \ \ \ \ +2\left(\delta_{\alpha+\beta,0}
-q^{2\alpha}\delta_{\alpha-\beta,0}\right)\frac{\partial}{\partial
z}x^{-1}\delta\left(\frac{z}{x}\right)\cdot {\bf c}.
\end{eqnarray}

Define a linear character $\chi_{q}^{D}: \Z_{2}\times \Z\rightarrow
\C^{\times}$ by
\begin{eqnarray}
\chi_{q}^{D} (\tau)=1\ \mbox{ and } \ \chi_{q}^{D}(\sigma_{r})=q^{r}\ \mbox{
for }r\in \Z.
\end{eqnarray}
Now, we have:

\bp{ptypeD} Lie algebra $\widehat{D}_{\hbar}$ is isomorphic to the
$(\Z_{2}\times \Z,\chi_{q}^{D})$-covariant algebra of affine Lie
algebra $\widehat{\A}$ with $c=\frac{1}{2}{\bf k}$ and
\begin{eqnarray}
D_{\alpha,m}=q^{\alpha}\overline{G_{\alpha,0}\otimes t^{m}} \ \ \
\mbox{ for }\alpha,m\in \Z.
\end{eqnarray}
\ep

\begin{proof} For $\alpha,m\in \Z$, set $D'_{\alpha,m}=q^{-\alpha}D_{\alpha,m}$.
Then we have
\begin{eqnarray}
D'_{-\alpha,m}=-D'_{\alpha,m}
\end{eqnarray}
and
\begin{eqnarray}\label{eD'-relation}
&&[D'_{\alpha,m},D'_{\beta,n}]=\left(q^{m\beta-n\alpha}-q^{n\alpha-m\beta}\right)D'_{\alpha+\beta,m+n}\nonumber\\
&&\ \ \ \ \ \
+\left(q^{m\beta+n\alpha}-q^{-m\beta-n\alpha}\right)D'_{\alpha-\beta,m+n}
+2m(\delta_{\alpha+\beta,0}-\delta_{\alpha-\beta,0})\delta_{m+n,0}\cdot
{\bf c}\ \hspace{1cm}
\end{eqnarray}
for $\alpha,\beta,m,n\in \Z$. We see that $\widehat{D}_{\hbar}$ can
be defined alternatively as the Lie algebra with generators $D'_{\alpha,m}$ for
$\alpha,m\in \Z$, subject to these relations.

On the other hand, consider the covariant algebra
$\widehat{\A}[\Z_{2}\times \Z]$ with respect to the linear character
$\chi_{q}^{D}$.  Let $\alpha,\beta, r,m,n\in \Z$. Just as with
$\hat{A}_{\hbar}$, we have
\begin{eqnarray*}
\sum_{r\in
\Z}\chi_{q}^{D}(\sigma_{r})^{m}\overline{[\sigma_{r}G_{\alpha,0},G_{\beta,0}]\otimes
t^{m+n}}
=\left(q^{m\beta-n\alpha}-q^{n\alpha-m\beta}\right)\overline{G_{\alpha+\beta,0}\otimes
t^{m+n}}.
\end{eqnarray*}
As $\chi_{q}^{D}(\tau\sigma_{r})=q^{r}$ for $r\in \Z$, we have
\begin{eqnarray*}
&&\sum_{r\in
\Z}\chi_{q}^{D}(\tau\sigma_{r})^{m}\overline{[\tau\sigma_{r}G_{\alpha,0},G_{\beta,0}]\otimes
t^{m+n}}\\
&=&-\sum_{r\in \Z}q^{rm}\overline{[G_{-\alpha,r},G_{\beta,0}]\otimes t^{m+n}}\\
&=&-\sum_{r\in
\Z}q^{rm}\left(\delta_{r,\beta-\alpha}\overline{G_{\beta-\alpha,-\alpha}\otimes
t^{m+n}}
-\delta_{r,\alpha-\beta}\overline{G_{\beta-\alpha,\alpha}\otimes
t^{m+n}}\right)\\
&=&\sum_{r\in
\Z}q^{rm}\left(\delta_{r,\beta-\alpha}\overline{\tau\sigma_{-\alpha}G_{\alpha-\beta,0}\otimes
t^{m+n}}-
\delta_{r,\alpha-\beta}\overline{\tau\sigma_{\alpha}G_{\alpha-\beta,0}\otimes
t^{m+n}}\right)\\
&=&\sum_{r\in
\Z}\left(\delta_{r,\beta-\alpha}q^{rm}q^{\alpha(m+n)}\overline{G_{\alpha-\beta,0}\otimes
t^{m+n}}-
\delta_{r,\alpha-\beta}q^{rm}q^{-\alpha(m+n)}\overline{G_{\alpha-\beta,0}\otimes
t^{m+n}}\right)\\
&=&\left(q^{m\beta+n\alpha}-q^{-m\beta-n\alpha}\right)\overline{G_{\alpha-\beta,0}\otimes
t^{m+n}}.
\end{eqnarray*}
We also have
\begin{eqnarray*}
&&\chi_{q}^{D}(\sigma_{r})^{m}\<\sigma_{r}G_{\alpha,0},G_{\beta,0}\>=q^{rm}\<G_{\alpha,r},G_{\beta,0}\>
=q^{rm}\delta_{\alpha+\beta,0}\delta_{r,0}=\delta_{\alpha+\beta,0}\delta_{r,0},\\
&&\chi_{q}^{D}(\tau\sigma_{r})^{m}\<\tau\sigma_{r}G_{\alpha,0},G_{\beta,0}\>=-q^{rm}\<G_{-\alpha,r},G_{\beta,0}\>
=-q^{rm}\delta_{\beta-\alpha,0}\delta_{r,0}=-\delta_{\alpha-\beta,0}\delta_{r,0}.
\end{eqnarray*}
Using all these relations we get
\begin{eqnarray}
&&[\overline{G_{\alpha,0}\otimes t^{m}},\overline{G_{\beta,0}\otimes
t^{n}}]=\left(q^{m\beta-n\alpha}-q^{n\alpha-m\beta}\right)\overline{G_{\alpha+\beta,0}\otimes
t^{m+n}}\nonumber\\
&&\ \ \ \
+\left(q^{m\beta+n\alpha}-q^{-m\beta-n\alpha}\right)\overline{G_{\alpha-\beta,0}\otimes
t^{m+n}}+m\left(\delta_{\alpha+\beta,0}-\delta_{\alpha-\beta,0}\right)\delta_{m+n,0}
{\bf k}.\ \ \ \ \ \
\end{eqnarray}
Furthermore, for $\alpha,m\in \Z$ we have
\begin{eqnarray*}
\overline{G_{-\alpha,0}\otimes t^{m}}=-\overline{\tau
G_{\alpha,0}\otimes
t^{m}}=-\chi_{q}^{D}(\tau)^{-n}\overline{G_{\alpha,0}\otimes t^{m}}
=-\overline{G_{\alpha,0}\otimes t^{m}}.
\end{eqnarray*}
Now, the assertion follows.
\end{proof}

Recall  from Definition \ref{dtau} the order-$2$ automorphism $\tau$ of $\gl_{\infty}$, which preserves $\A$.
Denote by $\A^{\tau}$ the Lie subalgebra of $\tau$-fixed points in
$\A$:
\begin{eqnarray}
\A^{\tau}=\{ a\in \A\   |\   \tau(a)=a\}.
\end{eqnarray}
For $\alpha,m\in \Z$, set
\begin{eqnarray}
G_{\alpha,m}^{\tau}=G_{\alpha,m}-G_{-\alpha,m}\in \A^{\tau}.
\end{eqnarray}
Note that  $G_{\alpha,m}^{\tau}$
for $\alpha\ge 1,\ m\in \Z$ form a basis of $\A^{\tau}$. Using (\ref{eGalpham}) we get
\begin{eqnarray}
[G_{\alpha,m}^{\tau},G_{\beta,n}^{\tau}]&=&\delta_{m-\alpha,n+\beta}G_{\alpha+\beta,n+\alpha}^{\tau}
-\delta_{m+\alpha,n-\beta}G_{\alpha+\beta,\beta+m}^{\tau}\nonumber\\
&&+\delta_{m+\alpha,n+\beta}G_{\alpha-\beta,m-\beta}^{\tau}
-\delta_{m-\alpha,n-\beta}G_{\alpha-\beta,\alpha+n}^{\tau},
\end{eqnarray}
\begin{eqnarray}
\<G_{\alpha,m}^{\tau},G_{\beta,n}^{\tau}\>
=2\left(\delta_{\alpha+\beta,0}-\delta_{\alpha-\beta,0}\right)\delta_{m,n}
\end{eqnarray}
for $\alpha,\beta,m,n\in \Z$. We have the following variation of
Proposition \ref{ptypeD}:

\bp{pD-type} Lie algebra $\widehat{D}_{\hbar}$ is isomorphic to
the covariant algebra $\widehat{\A^{\tau}}[\Z]$ with
$$D_{\alpha,m}=q^{\alpha}\overline{G_{\alpha,0}^{\tau}\otimes t^{m}}
\ \ \mbox{ for }\alpha,m\in \Z$$  and with ${\bf c}={\bf k}$. \ep

\begin{proof} It is similar to the proof of Proposition \ref{ptypeD}.
Notice that ${\bf c}$ and $D'_{\alpha,m}$ for $\alpha\ge 1,\ m\in
\Z$ form a basis of $\widehat{D}_{\hbar}$, where
$D'_{\alpha,m}=q^{-\alpha}D_{\alpha,m}$.  On the other hand, ${\bf
k}$ and $\overline{G_{\alpha,0}^{\tau}\otimes t^{m}}$ for $\alpha\ge
1,\ m\in \Z$ form a basis of $\widehat{\A^{\tau}}[\Z]$. Then we have
a linear isomorphism $\theta: \widehat{D}_{\hbar}\rightarrow
\widehat{\A^{\tau}}[\Z]$ such that $\theta({\bf c})={\bf k}$ and
$$\theta(D_{\alpha,m})=q^{\alpha}\overline{\left(G_{\alpha,0}^{\tau}\otimes t^{m}\right)}\ \ \
\mbox{ for }\alpha,m\in \Z,$$ noticing that
$D_{-\alpha,m}'=-D_{\alpha,m}'$ in $\widehat{D}_{\hbar}$ and
$\overline{G_{-\alpha,0}^{\tau}\otimes
t^{n}}=-\overline{G_{\alpha,0}^{\tau}\otimes t^{n}}$ in
$\widehat{\A^{\tau}}[\Z]$. Let $\alpha,\beta,m,n\in \Z$. As
$\sigma_{r}(G_{\alpha,0}^{\tau})=G_{\alpha,r}^{\tau}$ and
$\chi_{q}(\sigma_{r})=q^{r}$ for $r\in \Z$, from Proposition
\ref{ptwisted-Lie} we have
\begin{eqnarray}
&&\left[G_{\alpha,0}^{\tau}\otimes t^{m}, G_{\beta,0}^{\tau}\otimes
t^{n}\right]_{\Z}\nonumber\\
&=&\sum_{r\in
\Z}q^{rm}\left([G_{\alpha,r}^{\tau},G_{\beta,0}^{\tau}]\otimes
t^{m+n}+m\delta_{m+n,0}\<G_{\alpha,r}^{\tau},G_{\beta,0}^{\tau}\>{\bf
k}\right).
\end{eqnarray}
Furthermore, we have
\begin{eqnarray*}
&&\sum_{r\in
\Z}q^{rm}\overline{[G_{\alpha,r}^{\tau},G_{\beta,0}^{\tau}]\otimes t^{m+n}}\\
&=&\sum_{r\in
\Z}q^{rm}\overline{\left(\delta_{r,\alpha+\beta}G_{\alpha+\beta,\alpha}^{\tau}
-\delta_{r,-\alpha-\beta}G_{\alpha+\beta,-\alpha}^{\tau}+\delta_{r,\beta-\alpha}G_{\alpha-\beta,-\alpha}^{\tau}
-\delta_{r,\alpha-\beta}G_{\alpha-\beta,\alpha}^{\tau}\right)
\otimes t^{m+n}}\\
&=&\overline{\left(q^{m(\alpha+\beta)}G_{\alpha+\beta,\alpha}^{\tau}
-q^{-m(\alpha+\beta)}G_{\alpha+\beta,-\alpha}^{\tau}+q^{m(\beta-\alpha)}G_{\alpha-\beta,-\alpha}^{\tau}
-q^{m(\alpha-\beta)}G_{\alpha-\beta,\alpha}^{\tau}\right)
\otimes t^{m+n}}\\
&=&\overline{\left(q^{m(\alpha+\beta)}\sigma_{\alpha}G_{\alpha+\beta,0}^{\tau}
-q^{-m(\alpha+\beta)}\sigma_{-\alpha}G_{\alpha+\beta,0}^{\tau}\right)\otimes t^{m+n}}\\
&&+\overline{\left(q^{m(\beta-\alpha)}\sigma_{-\alpha}G_{\alpha-\beta,0}^{\tau}
-q^{m(\alpha-\beta)}\sigma_{\alpha}G_{\alpha-\beta,0}^{\tau}\right)
\otimes t^{m+n}}\\
&=&\left(q^{m\beta-n\alpha}-q^{n\alpha-m\beta}\right)\overline{G_{\alpha+\beta,0}^{\tau}\otimes
t^{m+n}}+\left(q^{m\beta+n\alpha}-q^{-m\beta-n\alpha}\right)\overline{G_{\alpha-\beta,0}\otimes
t^{m+n}}
\end{eqnarray*}
and
\begin{eqnarray*}
\sum_{r\in \Z}\< G_{\alpha,r}^{\tau}, G_{\beta,0}^{\tau}\>
=\sum_{r\in
\Z}2(\delta_{\alpha+\beta,0}-\delta_{\alpha-\beta,0})\delta_{r,0}
=2(\delta_{\alpha+\beta,0}-\delta_{\alpha-\beta,0}).
\end{eqnarray*}
Then we obtain
\begin{eqnarray*}
&&\left[\overline{G_{\alpha,0}^{\tau}\otimes t^{m}},
\overline{G_{\beta,0}^{\tau}\otimes t^{n}}\right]\nonumber\\
&=&\left(q^{m\beta-n\alpha}-q^{n\alpha-m\beta}\right)\overline{G_{\alpha+\beta,0}^{\tau}\otimes
t^{m+n}}+\left(q^{m\beta+n\alpha}-q^{-m\beta-n\alpha}\right)\overline{G_{\alpha-\beta,0}\otimes
t^{m+n}}\\
&&+2m\delta_{m+n,0}(\delta_{\alpha+\beta,0}-\delta_{\alpha-\beta,0}){\bf
k}.
\end{eqnarray*}
It now follows from this and (\ref{eD'-relation}) that $\theta$ is a
Lie algebra isomorphism.
\end{proof}

\section{Interplay between modules for trigonometric Lie algebras and quasi modules for vertex algebras}
In this section, we first recall from \cite{li-gamma} (cf. \cite{li-tlie})
the notion of
vertex $\Gamma$-algebra and the notion of equivariant quasi
module for a vertex $\Gamma$-algebra, and then as the main results, we show
that restricted modules of level $\ell$ for trigonometric Lie
algebras are equivariant quasi modules for the vertex algebra
$V_{\widehat{\A}}(\ell,0)$ with $\Gamma$ equal to the automorphism group $\Z$
or $\Z_{2}\times \Z$.

We begin by recalling from \cite{li-gamma} the notion of quasi module for
a vertex algebra.

\bd{dquasi-module-1} {\em Let $V$ be a vertex algebra. A {\em quasi
$V$-module} is a vector space $W$ equipped with a linear map
\begin{eqnarray*}
Y_{W}(\cdot,x):&& V\rightarrow \Hom (W,W((x)))\subset (\End
W)[[x,x^{-1}]]\\
&&v\mapsto Y_{W}(v,x),
\end{eqnarray*}
satisfying the conditions that $Y_{W}({\bf 1},x)=1_{W}$ (the
identity operator on $W$) and that for $u,v\in V$, there exists a
nonzero polynomial $p(x_{1},x_{2})$ such that
\begin{eqnarray}\label{epjacobi}
&&x_{0}^{-1}\delta\left(\frac{x_{1}-x_{2}}{x_{0}}\right)p(x_{1},x_{2})Y_{W}(u,x_{1})Y_{W}(v,x_{2})
\nonumber\\
&&\hspace{1cm}
-x_{0}^{-1}\delta\left(\frac{x_{2}-x_{1}}{-x_{0}}\right)p(x_{1},x_{2})Y_{W}(v,x_{2})Y_{W}(u,x_{1})\nonumber\\
&=&x_{2}^{-1}\delta\left(\frac{x_{1}-x_{0}}{x_{2}}\right)p(x_{1},x_{2})Y_{W}(Y(u,x_{0})v,x_{2})
\end{eqnarray}
(the {\em quasi Jacobi identity}).}
\ed

It is clear that the notion of quasi module generalizes that of module.  In
the definition of a module, the Jacobi identity can be replaced by
weak associativity (cf. \cite{li-local}, \cite{ll}). The following
is an analog for quasi modules, whose proof is straightforward (cf. \cite{ltw-twisted}):

\bl{lquasi-assoc} Let $V$ be a vertex algebra. In the definition of
a quasi $V$-module, the quasi Jacobi identity can be equivalently
replaced with the property that for $u,v\in V$, there exists a
nonzero polynomial $q(x_{1},x_{2})$ such that
\begin{eqnarray}\label{eqywuv}
q(x_{1},x_{2})Y_{W}(u,x_{1})Y_{W}(v,x_{2})\in \Hom
\left(W,W((x_{1},x_{2}))\right)
\end{eqnarray}
and
\begin{eqnarray}\label{eqywuv-assoc}
q(x_{2}+x_{0},x_{2})Y_{W}(Y(u,x_{0})v,x_{2})=
\left(q(x_{1},x_{2})Y_{W}(u,x_{1})Y_{W}(v,x_{2})\right)|_{x_{1}=x_{2}+x_{0}}.
\end{eqnarray}
\el

The following is a simple fact that we frequently use:

\bl{lsimple-fact} Let $V$ be a vertex algebra and $(W,Y_{W})$ a quasi $V$-module.
For $u,v\in V$, if
$$h(x_{1},x_{2})Y_{W}(u,x_{1})Y_{W}(v,x_{2})\in \Hom (W,W((x_{1},x_{2})))$$
for some series $h(x_{1},x_{2})\in \C((x_{1},x_{2}))$, then
\begin{eqnarray}\label{ehywuv}
h(x_{2}+x_{0},x_{2})Y_{W}(Y(u,x_{0})v,x_{2})=\left(h(x_{1},x_{2})Y_{W}(u,x_{1})Y_{W}(v,x_{2})\right)|_{x_{1}=x_{2}+x_{0}}.
\end{eqnarray}
\el

\begin{proof} Let $q(x_{1},x_{2})$ be a nonzero polynomial such that
(\ref{eqywuv}) and (\ref{eqywuv-assoc}) hold. Then
\begin{eqnarray*}
&&h(x_{2}+x_{0},x_{2})q(x_{2}+x_{0},x_{2})Y_{W}(Y(u,x_{0})v,x_{2})\nonumber\\
&=&\left(h(x_{1},x_{2})q(x_{1},x_{2})Y_{W}(u,x_{1})Y_{W}(v,x_{2})\right)|_{x_{1}=x_{2}+x_{0}}\nonumber\\
&=&q(x_{2}+x_{0},x_{2})\left(h(x_{1},x_{2})Y_{W}(u,x_{1})Y_{W}(v,x_{2})\right)|_{x_{1}=x_{2}+x_{0}}.
\end{eqnarray*}
Notice that we are allowed to cancel the factor $q(x_{2}+x_{0},x_{0})$
(a nonzero element of $\C((x_{2}))((x_{0}))$),
and by doing so we obtain (\ref{ehywuv}).
\end{proof}

The following notion of vertex $\Gamma$-algebra was introduced in \cite{li-gamma}:

\bd{dgamma-va} {\em Let $\Gamma$ be a group. A {\em vertex
$\Gamma$-algebra} is a vertex algebra $V$ equipped with two group
homomorphisms
\begin{eqnarray*}
&&R: \Gamma \rightarrow {\rm GL}(V);\ \  g\mapsto R_{g}\\
&&\phi: \Gamma \rightarrow \C^{\times},
\end{eqnarray*}
satisfying the condition that $R_{g}({\bf 1})={\bf 1}$,
\begin{eqnarray}\label{egenerating-vgamma-algebra}
R_{g}Y(v,x)=Y(R_{g}(v),\phi(g)^{-1}x)R_{g} \ \ \mbox{ for }g\in
\Gamma,\ v\in V.
\end{eqnarray} } \ed

\br{rxaZ-fraded} {\em Let $V$ be a $\Z$-graded vertex algebra in
the sense that $V$ is a vertex algebra equipped with a $\Z$-grading
 $V=\oplus_{n\in \Z}V_{(n)}$ such that ${\bf 1}\in V_{(0)}$ and
\begin{eqnarray}
u_{m}V_{(n)}\subset V_{(n+k-m-1)}\ \ \mbox{ for }u\in V_{(k)},\ m,n,
k\in \Z.
\end{eqnarray}
Define a linear operator $L(0)$ on $V$ by $L(0)v=nv$ for $v\in
V_{(n)}$ with $n\in \Z$. Suppose that $\Gamma$ is a group of
automorphisms of $V$ that preserve the $\Z$-grading. Let $\phi:
\Gamma\rightarrow \C^{\times}$ be any linear character. For $g\in
\Gamma$, set $R_{g}=\phi(g)^{-L(0)}g$. Then $V$ becomes a vertex
$\Gamma$-algebra (see \cite{li-gamma}).} \er

\bd{dquasi-module-2} {\em Let $V$ be a vertex $\Gamma$-algebra. An
{\em equivariant quasi $V$-module} is a quasi module
$(W,Y_{W})$ for $V$ viewed as a vertex algebra, satisfying the
conditions that
\begin{eqnarray}
Y_{W}(R_{g}v,x)=Y_{W}(v,\phi(g)x)\ \ \ \mbox{ for }g\in \Gamma,\
v\in V
\end{eqnarray}
and  that for $u,v\in V$,
(\ref{epjacobi}) holds with some polynomial $p(x_{1},x_{2})$ of the form
\begin{eqnarray}
p(x_{1},x_{2})=(x_{1}-\phi(g_{1})x_{2})\cdots
(x_{1}-\phi(g_{k})x_{2})
\end{eqnarray}
for some (not necessarily distinct) $g_{1},\dots,g_{k}\in \Gamma$.} \ed

To emphasize the dependence on the group $\Gamma$ and the linear character $\chi$,
we shall also use the term $(\Gamma,\chi)$-equivariant quasi $V$-module.

The following result was obtained in \cite{li-tlie} (Proposition
2.13):

\bp{pcommutator-formula} Let $V$ be a vertex $\Gamma$-algebra, let $\psi: \phi(\Gamma)\rightarrow \Gamma$
be a section of the linear character $\phi: \Gamma\rightarrow \C^{\times}$, and
let $(W,Y_{W})$ be an equivariant quasi $V$-module. For
$u,v\in V$, we have
\begin{eqnarray}
&&[Y_{W}(u,x_{1}),Y_{W}(v,x_{2})]\nonumber\\
&=&\sum_{\alpha\in
\phi(\Gamma)}\Res_{x_{0}}x_{1}^{-1}\delta\left(\frac{\alpha
x_{2}+x_{0}}{x_{1}}\right)Y_{W}(Y(R_{\psi(\alpha)}u,\alpha^{-1}x_{0})v,x_{2}),
\end{eqnarray}
which is a finite sum. \ep

As an immediate consequence we have:

\bc{cqm-wcomm}
Under all the assumptions in Proposition \ref{pcommutator-formula}, let
$u,v\in V$. Set
$$S_{u,v}=\{ g\in {\rm Im} (\psi) \subset \Gamma\ |\ (gu)_{n}v\ne 0\ \ \mbox{  for some }n\ge 0\}.$$
For $g\in S_{u,v}$, let $k_{g}\in \N$ be such that $(gu)_{m}v=0$ for
$m\ge k_{g}$. Then  there are finitely many $g_{1},\dots,g_{r}\in S_{u,v}$ such that
\begin{eqnarray}
\left(\prod_{i=1}^{r}
(x_{1}-\phi(g_{i})x_{2})^{k_{g_{i}}}\right)[Y_{W}(u,x_{1}),Y_{W}(v,x_{2})]=0.
\end{eqnarray}
\ec

Now, consider a (general) affine Lie algebra $\hat{\g}$. For $a\in
\g$, form a generating function
$$a(x)=\sum_{n\in \Z}a(n)x^{-n-1}\in \hat{\g}[[x,x^{-1}]],$$
where $a(n)$ stands for $a\otimes t^{n}$.
The defining relations (\ref{eaffine-defining-relation-1}) can be
written as
\begin{eqnarray}
[a(x_{1}),b(x_{2})]=[a,b](x_{2})x_{1}^{-1}\delta\left(\frac{x_{2}}{x_{1}}\right)
+\<a,b\>\frac{\partial}{\partial
x_{2}}x_{1}^{-1}\delta\left(\frac{x_{2}}{x_{1}}\right){\bf k}
\end{eqnarray}
for $a,b\in \g$.
Let $\ell$ be a complex number. View $\C$ as a module for
subalgebra $\g\otimes \C[t]+\C {\bf k}$ with $\g\otimes \C[t]$
acting trivially and with ${\bf k}$ acting as scalar $\ell$. Then
form an induced module
\begin{eqnarray}
V_{\hat{\g}}(\ell,0)=U(\hat{\g})\otimes_{U(\g\otimes \C[t]+\C {\bf
k})} \C.
\end{eqnarray}
 Set ${\bf 1}=1\otimes 1\in
V_{\hat{\g}}(\ell,0)$ and identify $\g$ as a subspace of
$V_{\hat{\g}}(\ell,0)$ through linear map
 $$\g\ni a\mapsto
a(-1){\bf 1}\in V_{\hat{\g}}(\ell,0).$$ Then there
exists a vertex algebra structure on $V_{\hat{\g}}(\ell,0)$,
which is uniquely determined by the condition that
${\bf 1}$ is the vacuum vector and $Y(a,x)=a(x)$ for $a\in \g$.
Furthermore, $V_{\hat{\g}}(\ell,0)$ is a $\Z$-graded vertex algebra such that
$$V_{\hat{\g}}(\ell,0)_{(n)}=0\ \mbox{ for }n<0, \ \
V_{\hat{\g}}(\ell,0)_{(0)}=\C {\bf 1}, \ \mbox{ and }\
V_{\hat{\g}}(\ell,0)_{(1)}=\g.$$

Let $\Gamma$ be an automorphism group of $(\g, \<\cdot,\cdot\>)$.
A simple fact is that every automorphism of $\g$, which reserves
$\<\cdot,\cdot\>$, extends uniquely to an automorphism of vertex
algebra $V_{\hat{\g}}(\ell,0)$. In view of this,
we can and we should view $\Gamma$ as an automorphism group of
$V_{\hat{\g}}(\ell,0)$ as a $\Z$-graded vertex algebra. From Remark \ref{rxaZ-fraded}, for any given linear character $\phi$ of $\Gamma$, $V_{\hat{\g}}(\ell,0)$ has a canonical vertex $\Gamma$-algebra
structure.

Recall from Section 2 the affine Lie algebra
$\widehat{\A}$, where $\A$ is a subalgebra of $\gl_{\infty}$ with a
basis $\{ E_{m,n}\ |\ m+n\in 2\Z\}$. We also recall that
$$G_{\alpha,m}=E_{m+\alpha,m-\alpha}\ \ \mbox{ for }\alpha,m\in \Z.$$
For any complex number $\ell$, we have a $\Z$-graded vertex algebra
$V_{\widehat{\A}}(\ell,0)$ whose
degree-$1$ subspace is canonically identified with $\A$. In particular, we have
$$G_{\alpha,m}\in V_{\widehat{\A}}(\ell,0)_{(1)}\ \ \ \mbox{
for }\alpha,m\in \Z.$$ Every automorphism of $\A$ viewed as a Lie
algebra, which preserves $\<\cdot,\cdot\>$, extends canonically to
an automorphism of $\widehat{\A}$ and furthermore to an automorphism
of vertex algebra $V_{\widehat{\A}}(\ell,0)$, preserving the
$\Z$-grading. In view of this, the automorphisms $\sigma_{n}$ of
$\A$, which was defined in Definition \ref{dsigma-r}, and $\tau$
extend to automorphisms of $V_{\widehat{\A}}(\ell,0)$. Then we view
$\Z_{2}\times \Z$ as an automorphism group of
$V_{\widehat{\A}}(\ell,0)$, preserving the $\Z$-grading. For $n\in
\Z$, set
\begin{eqnarray}
R_{n}=q^{-nL(0)}\sigma_{n}\in {\rm GL}(V_{\widehat{\A}}(\ell,0)),
\end{eqnarray}
recalling that $L(0)$ denotes the degree operator of
$V_{\widehat{\A}}(\ell,0)$. Equipped with these structures,
$V_{\widehat{\A}}(\ell,0)$ becomes a vertex $\Z$-algebra.

We formulate the following routine notions:

\bd{drestricted} {\em An $\widehat{A}_{\hbar}$-module $W$ is said to be
of {\em level $\ell\in \C$} if $c$ acts on $W$ as scalar $\ell$, and
it is said to be {\em restricted} if for any $w\in W,\ \alpha\in
\Z$, $A_{\alpha,n}w=0$ for $n$ sufficiently large.} \ed

As the main result of this section we have:

\bt{ttypeA-quasimodule} Assume that $q$ is not a root of unity and
let $\ell\in \C$. Then for any restricted $\widehat{A}_{\hbar}$-module
$W$ of level $\ell$, there exists an equivariant quasi
$V_{\widehat{\A}}(\ell,0)$-module structure $Y_{W}(\cdot,x)$ on $W$,
which is uniquely determined by
\begin{eqnarray}
Y_{W}(G_{\alpha,m},x)=q^{m}A_{\alpha}(q^{m}x) \ \ \ \mbox{ for
}\alpha,m\in \Z.
\end{eqnarray}
On the other hand, for any equivariant quasi
$V_{\widehat{\A}}(\ell,0)$-module $(W,Y_{W})$, $W$ becomes a
restricted $\widehat{A}_{\hbar}$-module of level $\ell$ with
$$A_{\alpha}(z)=Y_{W}(G_{\alpha,0},z)\ \ \ \mbox{ for }\alpha\in
\Z.$$ \et

\begin{proof} Let $W$ be a restricted
$\widehat{A}_{\hbar}$-module of level $\ell$. In view of Proposition
\ref{ptypeA}, $W$ is a restricted $\widehat{\A}[\Z]$-module of level
$\ell$ with
$$\overline{G_{\alpha,0}\otimes t^{m}}=A_{\alpha,m}\ \ \mbox{ for
}\alpha,m\in \Z.$$ Note that as $q$ is not a root of unity,
$\chi_{q}:\Z\rightarrow \C^{\times}$ is one-to-one. By Theorem 4.9
of \cite{li-tlie}, there exists an equivariant quasi
$V_{\widehat{\A}}(\ell,0)$-module structure $Y_{W}(\cdot,x)$ on $W$,
which is uniquely determined by
\begin{eqnarray*}
Y_{W}(G_{\alpha,m},x)=\overline{G_{\alpha,m}(x)} \ \ \ \mbox{ for
}\alpha,m\in \Z.
\end{eqnarray*}
For $\alpha,m\in \Z$, we have
$$Y_{W}(G_{\alpha,m},x)=\overline{G_{\alpha,m}(x)}
=\overline{\sigma_{m}G_{\alpha,0}(x)}=q^{m}\overline{G_{\alpha,0}(q^{m}x)}=q^{m}A_{\alpha}(q^{m}x).$$
Consequently, there exists an equivariant quasi
$V_{\widehat{\A}}(\ell,0)$-module structure $Y_{W}(\cdot,x)$ on $W$
such that
\begin{eqnarray}
Y_{W}(G_{\alpha,m},x)=q^{m}A_{\alpha}(q^{m}x) \ \ \ \mbox{ for
}\alpha,m\in \Z.
\end{eqnarray}
The other direction follows from Proposition \ref{ptypeA} and
Theorem 4.9 of \cite{li-tlie}.
\end{proof}

Recall the linear character $\chi_{q}^{B}: \Z_{2}\times
\Z\rightarrow \C^{\times}$ which was defined in Section 2 by
$$\chi_{q}^{B}(\tau)=-1\ \ \mbox{ and }\ \ \chi_{q}^{B}(\sigma_{n})=q^{n}\ \ \mbox{
for }n\in \Z.$$
Note that if $q$ is not a root of unity, then $\chi_{q}$ is one-to-one.
View $V_{\widehat{\A}}(2\ell,0)$ as a vertex $(\Z_{2}\times \Z)$-algebra.
Combining Proposition \ref{ptypeB} with
Theorem 4.9 of \cite{li-tlie} we immediately have:

\bt{ttypeB-quasi-module} Assume that $q$ is not a root of unity and
let $\ell\in \C$. Then for any restricted $\widehat{B}_{\hbar}$-module
$W$ of level $\ell$, there exists a $(\Z_{2}\times
\Z,\chi_{q}^{B})$-equivariant quasi
$V_{\widehat{\A}}(2\ell,0)$-module structure $Y_{W}(\cdot,x)$ on
$W$, which is uniquely determined by
\begin{eqnarray}
Y_{W}(G_{\alpha,m},x)=q^{m}B_{\alpha}(q^{m}x) \ \ \ \mbox{ for
}\alpha,m\in \Z.
\end{eqnarray}
On the other hand, every $(\Z_{2}\times
\Z,\chi_{q}^{B})$-equivariant quasi
$V_{\widehat{\A}}(2\ell,0)$-module $W$ is a restricted $\widehat{B}_{\hbar}$-module of level $\ell$ with
$$B_{\alpha}(z)=Y_{W}(G_{\alpha,0},z)\ \ \ \mbox{ for }\alpha\in
\Z.$$
\et

Recall that $\tau(E_{m,n})=-E_{n,m}$ for $m,n\in \Z$, where $\tau$
is an order $2$ automorphism of the Lie algebra $\gl_{\infty}$,
which preserves the bilinear form. Furthermore, $\tau$ preserves
$\A$; $\tau(G_{\alpha,m})=-G_{-\alpha,m}$ for $\alpha,m\in \Z$.
Recall that $\A^{\tau}$ denotes the Lie subalgebra of $\tau$-fixed
points in $\A$. Fix a linear character $\chi_{q}^{D}: \Z\rightarrow
\C^{\times}$ defined by
$$\chi_{q}^{D}(\sigma_{n})=q^{n}\ \ \
\mbox{ for }n\in \Z.$$

We have:

\bt{ttypeD-quasi-module} Assume that $q$ is not a root of unity and
let $\ell\in \C$. Then for any restricted $\widehat{D}_{\hbar}$-module
$W$ of level $\ell$, there exists a $(\Z,\chi_{q}^{D})$-equivariant
quasi $V_{\widehat{\A^{\tau}}}(\ell,0)$-module structure
$Y_{W}(\cdot,x)$ on $W$,  which is uniquely determined by
\begin{eqnarray}
Y_{W}(G_{\alpha,m}^{\tau},x)=q^{-\alpha+m}D_{\alpha}(q^{m}x) \ \ \
\mbox{ for }\alpha,m\in \Z,
\end{eqnarray}
where $G_{\alpha,m}^{\tau}=G_{\alpha,m}-G_{-\alpha,m}\in \A^{\tau}$.
On the other hand, given a $(\Z,\chi_{q}^{D})$-equivariant quasi
$V_{\widehat{\A^{\tau}}}(\ell,0)$-module $(W,Y_{W})$, one has a
restricted $\widehat{D}_{\hbar}$-module structure of level $\ell$ on $W$
such that $D_{\alpha}(x)=q^{\alpha}Y_{W}(G_{\alpha,0}^{\tau},x)$ for
$\alpha\in \Z$. \et

\section{Vertex algebra $L_{\widehat{\A}}(\ell,0)$ and its modules}

In this section, we study the vertex algebra
$L_{\widehat{\A}}(\ell,0)$ with $\ell$ a positive integer and give a
characterization of $L_{\widehat{\A}}(\ell,0)$-modules.

We start with the $\Z$-graded vertex algebra $V_{\hat{\g}}(\ell,0)$
associated to an affine Lie algebra $\hat{\g}$ of a (possibly infinite dimensional) Lie algebra $\g$, where
$$V_{\hat{\g}}(\ell,0)_{(n)}=0\ \mbox{ for }n<0, \ \
V_{\hat{\g}}(\ell,0)_{(0)}=\C {\bf 1}, \ \mbox{ and }\
V_{\hat{\g}}(\ell,0)_{(1)}=\g.$$
Define a linear operator $\D$ on $V_{\hat{\g}}(\ell,0)$ by $\D v=v_{-2}{\bf 1}$
for $v\in V_{\hat{\g}}(\ell,0)$.
We have
\begin{eqnarray}\label{eDgrading}
\D V_{\hat{\g}}(\ell,0)_{(n)}\subset V_{\hat{\g}}(\ell,0)_{(n+1)}\ \ \ \mbox{ for }n\in \Z.
\end{eqnarray}

\bl{ldef-J} Let $J_{\hat{\g}}(\ell,0)$ be the (unique) maximal
graded $\hat{\g}$-submodule of $V_{\hat{\g}}(\ell,0)$. Then
$J_{\hat{\g}}(\ell,0)$ is a two-sided ideal of $V_{\hat{\g}}(\ell,0)$.
\el

\begin{proof} Since $\g$ generates
$V_{\hat{\g}}(\ell,0)$ as a vertex algebra, it follows that
$J_{\hat{\g}}(\ell,0)$ is a left ideal. To prove $J_{\hat{\g}}(\ell,0)$ is also a right ideal,
from \cite{ll} we must prove $\D J_{\hat{\g}}(\ell,0)\subset J_{\hat{\g}}(\ell,0)$.
Since $[\D,a_{n}]=-na_{n-1}$ for $a\in \g,\ n\in \Z$, it can be readily seen that
$J_{\hat{\g}}(\ell,0)+\D(J_{\hat{\g}}(\ell,0))$ is a graded $\hat{\g}$-submodule of $V_{\hat{\g}}(\ell,0)$. From (\ref{eDgrading}) the degree-zero homogeneous subspace of $J_{\hat{\g}}(\ell,0)+\D(J_{\hat{\g}}(\ell,0))$ is trivial.
Thus we have $J_{\hat{\g}}(\ell,0)+\D J_{\hat{\g}}(\ell,0)\subset J_{\hat{\g}}(\ell,0)$, proving
$\D J_{\hat{\g}}(\ell,0)\subset J_{\hat{\g}}(\ell,0)$. Therefore,
$J_{\hat{\g}}(\ell,0)$ is an ideal of $V_{\hat{\g}}(\ell,0)$.
\end{proof}

\bd{dL(l,0)} {\em Denote by $L_{\hat{\g}}(\ell,0)$ the quotient vertex algebra
of $V_{\hat{\g}}(\ell,0)$ modulo $J_{\hat{\g}}(\ell,0)$, which is a graded simple  vertex
algebra.}
\ed

Let $U$ be a $\g$-module.
We make $U$ a $(\g\otimes \C[t]+\C {\bf k})$-module by letting $\g\otimes t\C[t]$
act trivially and ${\bf k}$ act
as scalar $\ell$. Then form a generalized Verma $\hat{\g}$-module
$$M_{\hat{\g}}(\ell,U)=U(\hat{\g})\otimes _{U(\g\otimes \C[t]+\C {\bf k})}U,$$
which is naturally $\N$-graded. Define $L_{\hat{\g}}(\ell,U)$ to be the
quotient module of $M_{\hat{\g}}(\ell,U)$ modulo the maximal graded submodule
with trivial degree-zero subspace.

Now, specifying $\g$ to $\A$, where
$$\A={\rm span}\{ E_{m,n}\ |\ m+n\in 2\Z\}\subset \gl_{\infty},$$
we have $\Z$-graded vertex algebras
$V_{\widehat{\A}}(\ell,0)$ and $L_{\widehat{\A}}(\ell,0)$.
If $\ell=0$, we
see that $\A$ generates a proper graded ideal of
$V_{\widehat{\A}}(0,0)$. Consequently, $L_{\widehat{\A}}(0,0)=\C$.
If $\ell\ne 0$, from \cite{ll} we have
$L_{\widehat{\A}}(\ell,0)_{(1)}=\A$.

\br{rautomorphismL} {\em It is straightforward to show that
$J_{\widehat{\A}}(\ell,0)$ is stable under every automorphism of
$V_{\widehat{\A}}(\ell,0)$, which preserves the $\Z$-grading. Thus,
$J_{\widehat{\A}}(\ell,0)$ is stable under the action of group
$\Z_{2}\times \Z$  as defined in Section 3. Consequently,
$\Z_{2}\times \Z$ acts on $L_{\widehat{\A}}(\ell,0)$ by
automorphisms that preserve the $\Z$-grading.} \er

Notice that
\begin{eqnarray}
\A={\rm span}\{ E_{m,n}\ |\ m,n\in 2\Z\}\oplus {\rm span}\{
E_{p,q}\; |\ p,q\in 1+2\Z\},
\end{eqnarray}
which is isomorphic to $\gl_{\infty}\oplus \gl_{\infty}$ as an
associative algebra. As a Lie algebra, $\A$ has a Cartan subalgebra
\begin{eqnarray}
H={\rm span}\{ E_{m,m}\ |\ m\in \Z\}.
\end{eqnarray}
Recall
$$G_{\alpha,m}=E_{m+\alpha,m-\alpha}\in \A \ \ \mbox{ for }\alpha,m\in \Z.$$
Then $H={\rm span}\{G_{0,m}\ |\ m\in \Z\}$ and  for any $\alpha,m\in
\Z$ with $\alpha\ne 0$, $G_{\alpha,m}$
is a root vector of $\A$ viewed as a Lie algebra.

First, we have:

\bl{lsimple-va} For any complex number $\ell$,
$L_{\widehat{\A}}(\ell,0)$ is an irreducible $\widehat{\A}$-module
and a simple vertex algebra. \el

\begin{proof} Notice that the second assertion follows from the first one.
For the first assertion, we need to show that every
$\widehat{\A}$-submodule of $V_{\widehat{\A}}(\ell,0)$ is graded. If
$\ell=0$, we have $L_{\widehat{\A}}(\ell,0)=\C$, which is
irreducible. We now assume $\ell\ne 0$. Let $W$ be an
$\widehat{\A}$-submodule of $V_{\widehat{\A}}(\ell,0)$ and let $v\in
W$. Then there exists a finite interval $I$ of $\Z$ such that $v\in
U(\widehat{\A_{I}}){\bf 1}$.  where  $\A_{I}={\rm span}\{ E_{m,n}\ |\
m,n\in I,\ m+n\in 2\Z\}$.  We assume that $|I|$ is even.
Then $\A_{I}$ is a Lie subalgebra isomorphic to
the direct sum of $\sl_{\frac{1}{2}|I|}\oplus \sl_{\frac{1}{2}|I|}$ and a two-dimensional
abelian algebra. We choose $I$ such that $\ell\ne -\frac{1}{2}|I|$. It is clear
that $U(\widehat{\A_{I}}){\bf 1}$ is a graded vertex subalgebra of
$V_{\widehat{\A}}(\ell,0)$, generated by $\A_{I}$. In view of the
P-B-W theorem, we have $U(\widehat{\A_{I}}){\bf 1}\simeq
V_{\widehat{\A_{I}}}(\ell,0)$ as a $\Z$-graded
$\widehat{\A_{I}}$-module. As $\ell\ne 0$ and $\ell\ne -\frac{1}{2}|I|$ (the
negative dual Coxeter number of $\sl_{\frac{1}{2}|I|}$),
$V_{\widehat{\A_{I}}}(\ell,0)$ is a vertex operator algebra with the
standard Segal-Sugawara conformal vector $\omega_{I}$. Furthermore,
the $\Z$-grading of $V_{\widehat{\A_{I}}}(\ell,0)$ by the
$L_{I}(0)$-weight coincides with the $\Z$-grading of
$V_{\widehat{\A}}(\ell,0)$. As an $\widehat{\A_{I}}$-submodule, $W$
is stable under the action of $L_{I}(0)$. It follows that every
homogeneous component of $v$ lies in $W$.  This proves that $W$ is a
graded subspace. Therefore, $L_{\widehat{\A}}(\ell,0)$ is an
irreducible $\widehat{\A}$-module.
\end{proof}

\bd{dintegrable-module} {\em An $\widehat{\A}$-module $W$ is said to
be {\em integrable} if $H$ is semisimple on $W$ and if for any $\alpha,
m,k\in \Z$ with $\alpha\ne 0$, $G_{\alpha,m}(k)$ is locally
nilpotent on $W$.} \ed

\bl{lintegrability} Let $\ell$ be a complex number and let $W$ be a
nonzero restricted $\widehat{\A}$-module of level
$\ell$, on which $H$ is semisimple.
Then $W$ is integrable if and only if $\ell$ is a nonnegative integer and
for any $\alpha, m\in \Z$ with $\alpha\ne 0$,
$$G_{\alpha,m}(x)^{\ell+1}=0\ \ \mbox{ on
}\ W.$$  \el

\begin{proof} Assume that $W$ is integrable. Let $m,n\in \Z$ with $m+n\in 2\Z,\ m\ne n$.
We have a Lie algebra embedding of $\widehat{\sl(2,\C)}$ into
$\widehat{\A}$:
$${\bf k}\mapsto {\bf k},\ \
e\otimes t^{r}\mapsto E_{m,n}\otimes t^{r},\ \ f\otimes t^{r}\mapsto
E_{n,m}\otimes t^{r}, \ \ h\otimes t^{r}\mapsto
(E_{m,m}-E_{n,n})\otimes t^{r}$$ for $r\in \Z$. Then $W$ is
necessarily a restricted and integrable $\widehat{\sl(2,\C)}$-module
of level $\ell$ which must be a nonnegative integer. From
\cite{dlm}, $W$ is a direct sum of irreducible integrable highest
weight $\widehat{\sl(2,\C)}$-modules. It then follows from \cite{lp}
that
$$E_{m,n}(x)^{\ell+1}=E_{n,m}(x)^{\ell +1}=0\ \ \mbox{ on
}\ W.$$

The other direction essentially follows from \cite{dlm}: Assume that
$G_{\alpha,m}(x)^{\ell+1}=0$ on $W$ for $\alpha,m\in \Z$ with $\alpha\ne
0$. Let $w\in W$. We need to show that for $n\in \Z$,
$G_{\alpha,m}(n)^{q}w=0$ for some positive integer $q$. This will
follow from deduction on $n$ as follows. First, there exists an
integer $k$ such that $G_{\alpha,m}(p)w=0$ for $p>k$. Let $A$ be the
commutative subalgebra of $U(\widehat{\A})$, generated by
$G_{\alpha,m}(n)$ for $n\in \Z$, and set $E=Aw\subset W$. We see
that $G_{\alpha,m}(n)=0$ on $E$ for $n>k$. Assume that
$G_{\alpha,m}(n)$ for $n>k-r$ are all nilpotent on $E$ for some
nonnegative integer $r$. By extracting the coefficient of
$x^{-(\ell+1)(k-r+1)}$ from $G_{\alpha,m}(x)^{\ell+1}=0$ on $E$, we
get
$$G_{\alpha,m}(k-r)^{\ell+1}=G_{\alpha,m}(k-r)X\ \ \mbox{ on }E,$$
where $X$ is some polynomial of $G_{\alpha,m}(k-p)$ with $p<r$. It
follows that $G_{\alpha,m}(k-r)$ is also nilpotent on $E$.
This proves that $W$ is an integrable $\widehat{\A}$-module.
\end{proof}

For the vertex algebras $V_{\widehat{\A}}(\ell,0)$ and $L_{\widehat{\A}}(\ell,0)$, we have:

\bp{pquotient-L} Let $\ell$ be a nonnegative integer. Then the
maximal graded $\widehat{\A}$-submodule $J_{\widehat{\A}}(\ell,0)$
of $V_{\widehat{\A}}(\ell,0)$ is generated by
$G_{\alpha,m}(-1)^{\ell+1}{\bf 1}$ for $\alpha,m\in \Z$ with
$\alpha\ne 0$. Furthermore, $L_{\widehat{\A}}(\ell,0)$ is an
integrable $\widehat{\A}$-module.\ep

\begin{proof} Let $m,n\in \Z$ with $m+n\in 2\Z$ and $m\ne n$. Note that $E_{m,n}(-1)^{\ell+1}{\bf 1}$
is a homogeneous vector of degree $\ell+1$ in
$V_{\widehat{\A}}(\ell,0)$. We claim
\begin{eqnarray}
E_{p,q}(1)E_{m,n}(-1)^{\ell+1}{\bf 1}=0
\ \ \ \mbox{ for any }p,q\in \Z\ \mbox{ with }p+q\in 2\Z,\ p\ne q.
\end{eqnarray}
If $p\ne n$ and $q\ne m$, we have $[E_{p,q}(1),E_{m,n}(-1)]=0$, so that
$$E_{p,q}(1)E_{m,n}(-1)^{\ell+1}{\bf 1}=E_{m,n}(-1)^{\ell+1}E_{p,q}(1){\bf 1}=0.$$
Assume $(p,q)=(n,m)$. Notice that $E_{n,m}(1), E_{m,n}(-1),
E_{n,n}(0)-E_{m,m}(0)+\ell$ linearly span a Lie subalgebra
isomorphic to $\sl(2,\C)$ with
$$e=E_{m,n}(1),\ \ f=E_{n,m}(-1),\ \ h=E_{n,n}(0)-E_{m,m}(0)+\ell.$$
Since
$$E_{n,m}(1){\bf 1}=0,\ \ \
(E_{n,n}(0)-E_{m,m}(0)+\ell){\bf 1}=\ell {\bf 1},$$
we have $E_{n,m}(1)E_{m,n}(-1)^{\ell+1}{\bf 1}=0$.

Consider the case where $q= m$ and $p\ne n$. We have $[E_{p,q}(1),E_{m,n}(-1)]=E_{p,n}(0)$.
As $[E_{p,n}(0),E_{m,n}(-1)]=0$, we have
$$E_{p,n}(0)E_{m,n}(-1)^{k}{\bf 1}=E_{m,n}(-1)^{k}E_{p,n}(0){\bf 1}=0$$
for any nonnegative integer $k$. It follows from induction that
$$E_{p,q}(1)E_{m,n}(-1)^{k+1}{\bf 1}=0\ \ \ \mbox{
for all nonnegative integers }k.$$ For the case $p=n$ and
$q\ne m$, the proof is similar. This proves our claim. It then follows (from
the P-B-W theorem) that $U(\widehat{\A})E_{m,n}(-1)^{\ell+1}{\bf 1}$
is a proper graded submodule. Therefore
$$E_{m,n}(-1)^{\ell+1}{\bf 1}\in J_{\widehat{\A}}(\ell,0).$$

Denote by $J'$ the $\widehat{\A}$-submodule of
$V_{\widehat{\A}}(\ell,0)$, generated by $E_{m,n}(-1)^{\ell+1}{\bf
1}$ for $m,n\in \Z$ with $m+n\in 2\Z,\ m\ne n$. Set
$V=V_{\widehat{\A}}(\ell,0)/J'$. It follows that for $m,n\in \Z$
with $m+n\in 2\Z,\ m\ne n$, $E_{m,n}(-1)$ is locally nilpotent on
$V$. For $k\ge 0$, it is clear that $E_{m,n}(k)$ is locally
nilpotent on $V$. Let $I$ be any finite interval of $\Z$ and let
$\A_{I}$ be the corresponding subalgebra of $\A$. From \cite{kac1},
$V$ is an integrable $\widehat{\A_{I}}$-module. As $I$ is
arbitrary, it follows
that $V$ is an integrable $\widehat{\A}$-module. On the other hand,
 as a submodule of the integrable $\widehat{\A_{I}}$-module $V$,
 $U(\widehat{\A_{I}}){\bf 1}$ must be an irreducible
$\widehat{\A_{I}}$-module. It follows that $V$ is an irreducible
$\widehat{\A}$-module. Consequently, we have
$V=L_{\widehat{\A}}(\ell,0)$, proving $J'=J_{\widehat{\A}}(\ell,0)$.
Then the two assertions hold.
\end{proof}

An an immediate consequence of Lemma \ref{lintegrability} and
Proposition \ref{pquotient-L}, we have:

\bc{cva-integ}  Let $\ell$ be a complex number. Then
$L_{\widehat{\A}}(\ell,0)$ is an integrable $\widehat{\A}$-module if
and only if $\ell$ is a nonnegative integer. \ec

{\em For the rest of this section, we assume that $\ell$ is a
nonnegative integer.} Note that each $L_{\widehat{\A}}(\ell,0)$-module
is naturally a restricted $\widehat{\A}$-module of level $\ell$.
Furthermore, we have:

\bp{pmodule-classification-1} All $L_{\widehat{\A}}(\ell,0)$-modules are
exactly those restricted $\widehat{\A}$-modules $W$ of level $\ell$
such that
$$G_{\alpha,m}(x)^{\ell+1}=0\  \mbox{ on }W$$
for $\alpha,m\in \Z$ with $\alpha\ne 0$.
\ep

\begin{proof} Let $(W,Y_{W})$ be an $L_{\widehat{\A}}(\ell,0)$-module.
Then $W$ is a restricted $\widehat{\A}$-module of level $\ell$ with
$G_{\alpha,m}(x)=Y_{W}(G_{\alpha,m},x)$
for $\alpha,m\in \Z$ with $\alpha\ne 0$.
{}From Proposition \ref{pquotient-L} we have
$G_{\alpha,m}(-1)^{\ell+1}{\bf 1}=0$ in $L_{\widehat{\A}}(\ell,0)$.
By a result of Dong-Lepowsky in \cite{dl} we have
$$G_{\alpha,m}(x)^{\ell+1}=Y_{W}(G_{\alpha,m},x)^{\ell+1}=Y_{W}(G_{\alpha,m}(-1)^{\ell+1}{\bf
1},x)=0\  \mbox{ on }W.$$ Conversely, let $W$ be a restricted
$\widehat{\A}$-module level $\ell$, satisfying the very property.
First, $W$ is naturally a
$V_{\widehat{\A}}(\ell,0)$-module such that
$Y_{W}(G_{\alpha,m},x)=G_{\alpha,m}(x)$ for $\alpha,m\in \Z$.
Furthermore, for $\alpha\ne 0$ we have
$$Y_{W}\left(G_{\alpha,m}(-1)^{\ell+1}{\bf 1},x\right)=Y_{W}(G_{\alpha,m},x)^{\ell+1}=G_{\alpha,m}(x)^{\ell+1}=0.$$
Then it follows from Proposition \ref{pquotient-L} that $W$ is
naturally an $L_{\widehat{\A}}(\ell,0)$-module.
\end{proof}

As an immediate consequence of Lemma \ref{lintegrability}
and Proposition \ref{pmodule-classification-1}, we have:

\bc{cirreducible-modules} Let $W$ be a restricted
$\widehat{\A}$-module on which $H$ is semisimple. Then $W$ is an
$L_{\widehat{\A}}(\ell,0)$-module if and only if $W$ is an
integrable $\widehat{\A}$-module of level $\ell$. \ec

\br{raffine-voa}
{\em Let $\g$ be a finite-dimensional simple Lie algebra and let $\ell$ be a positive integer.
It was known (see \cite{fz}, \cite{dl}, \cite{li-local}) that
irreducible modules for $L_{\hat{\g}}(\ell,0)$ viewed as a vertex operator algebra  are exactly
integrable highest weight $\hat{\g}$-modules of level $\ell$.
Furthermore, it was proved in \cite{dlm} that
every  module for $L_{\hat{\g}}(\ell,0)$ viewed as a vertex algebra,
 is a direct sum of integrable highest weight $\hat{\g}$-modules of level $\ell$.}
\er

\section{Quasi modules for simple vertex algebra $L_{\widehat{\A}}(\ell,0)$ }

In this section, we determine irreducible $\Gamma_{q}$-equivariant quasi
$L_{\widehat{\A}}(\ell,0)$-modules
and we show that every
unitary quasifinite highest weight $\bar{\A}_{\hbar}$-module of level
$\ell$ is an irreducible  $\Gamma_{q}$-equivariant quasi $L_{\widehat{\A}}(\ell,0)$-module.

First, we establish a technical result.

\bl{lmodule-ab}
Let $\Gamma$ be a group equipped with a linear character $\phi: \Gamma\rightarrow \C^{\times}$.
Let $V$ be a vertex $\Gamma$-algebra, let $a,b\in V$ such that $a_{n}b=0$ for $n\ge 0$, and
let $(W,Y_{W})$ be a $\Gamma$-equivariant quasi $V$-module.
Then there exists a polynomial $q(z)$ with $q(1)=1$ such that
\begin{eqnarray}
q(x_{1}/x)[Y_{W}(a,x_{1}),Y_{W}(b,x)]=0.\label{eqab-locality}
\end{eqnarray}
Furthermore, for any such polynomial $q(z)$ we have
\begin{eqnarray}
Y_{W}(a_{-1}b,x) =
\left(q(x_{1}/x)Y_{W}(a,x_{1})Y_{W}(b,x)\right)|_{x_{1}=x}.\label{ea-1b}
\end{eqnarray}
\el

\begin{proof} Since $a_{n}b=0$ for $n\ge 0$, from Corollary \ref{cqm-wcomm}
there exist some  (possibly the same) nonzero and nonunit complex
numbers $\alpha_{1},\dots,\alpha_{k}$ such that
\begin{eqnarray*}
(x_{1}-\alpha_{1}x_{2})\cdots (x_{1}-\alpha_{k}x_{2})[Y_{W}(a,x_{1}),Y_{W}(b,x_{2})]=0.
\end{eqnarray*}
Set $q(z)=\prod_{i=1}^{k}\left(\frac{z-\alpha_{i}}{1-\alpha_{i}}\right)$. Then
$q(1)=1$ and (\ref{eqab-locality}) holds.
Note that (\ref{eqab-locality}) implies
$$q(x_{1}/x)Y_{W}(a,x_{1})Y_{W}(b,x)\in \Hom (W,W((x_{1},x))).$$
In view of Lemma \ref{lsimple-fact}, we have
\begin{eqnarray*}
q((x+x_{0})/x)Y_{W}(Y(a,x_{0})b,x)
= \left(q(x_{1}/x)Y_{W}(a,x_{1})Y_{W}(b,x)\right)|_{x_{1}=x+x_{0}}.
\end{eqnarray*}
As $Y(a,x_{0})b\in V[[x_{0}]]$, we can set $x_{0}=0$, and by doing so we obtain
(\ref{ea-1b}).
\end{proof}

Furthermore, we have:

\bp{pnilpotent} Let $V$ be a vertex $\Gamma$-algebra, let $a\in V$
such that $a_{n}a=0$ for $n\ge 0$, and let $(W,Y_{W})$ be any
$\Gamma$-equivariant quasi $V$-module. Then there exists a
polynomial $q(z)$ with $q(1)=1$ such that
\begin{eqnarray}\label{eqab-comm}
q(x_{1}/x)Y_{W}(a,x_{1})Y_{W}(a,x)=q(x_{1}/x)Y_{W}(a,x)Y_{W}(a,x_{1}).
\end{eqnarray}
Furthermore, for any such polynomial $q(z)$, we have
\begin{eqnarray}
&&Y_{W}\left((a_{-1})^{\ell+1}{\bf 1},x\right)\nonumber\\
&=&\left(P(x_{1},\dots, x_{\ell},x)Y_{W}(a,x_{1})\cdots
Y_{W}(a,x_{\ell})Y_{W}(a,x)\right)|_{x_{1}=x_{2}=\cdots =x_{\ell}=x},
\end{eqnarray}
where $\ell$ is a positive integer and
$P(x_{1},\dots,x_{\ell+1})=\prod_{1\le i<j\le
\ell+1}q(x_{i}/x_{j})$. \ep

\begin{proof}  Note that the first assertion and the second assertion with $\ell=1$
follow from Lemma \ref{lmodule-ab}. We now prove the second
assertion by induction on $\ell$. For $r\ge 2$, set
$$P_{r}(x_{1},\dots,x_{r})=\prod_{1\le i<j\le r}q(x_{i}/x_{j}).$$
Then
\begin{eqnarray}
&&P_{r}(x_{1},\dots,x_{r})Y_{W}(a,x_{\sigma(1)})Y_{W}(a,x_{\sigma(2)})\cdots
Y_{W}(a,x_{\sigma(r)})\nonumber\\
&=&
P_{r}(x_{1},\dots,x_{r})Y_{W}(a,x_{1})Y_{W}(a,x_{2})\cdots
Y_{W}(a,x_{r})
\end{eqnarray}
for any permutation $\sigma$ on $\{1,2,\dots,r\}$. {}From this we have
\begin{eqnarray*}
P_{r}(x_{1},\dots,x_{r})Y_{W}(a,x_{1})Y_{W}(a,x_{2})\cdots
Y_{W}(a,x_{r}) \in \Hom (W,W((x_{1},\dots,x_{r}))).
\end{eqnarray*}
In view of this, the substitution
\begin{eqnarray*}
\left(P_{r}(x_{1},\dots,x_{r})Y_{W}(a,x_{1})Y_{W}(a,x_{2})\cdots
Y_{W}(a,x_{r})\right)|_{x_{1}=\cdots =x_{r}=x}
\end{eqnarray*}
 exists in $\Hom (W,W((x)))$.

Now assume $\ell\ge 2$ and set $b=(a_{-1})^{\ell}{\bf 1}$.
Note that the assumption $a_{j}a=0$ for $j\ge 0$ is equivalent to that $[a_{m},a_{n}]=0$
for $m,n\in \Z$.
Then we have $a_{n}b=0$ for all $n\ge 0$ as $[a_{n},a_{-1}]=0$ and $a_{n}{\bf 1}=0$.
By induction hypothesis, we have
\begin{eqnarray}\label{e-hypothesis}
Y_{W}(b,x)
=\left(P_{\ell}(x_{2},\dots,x_{\ell},x)
Y_{W}(a,x_{2})\cdots Y_{W}(a,x_{\ell})Y_{W}(a,x)\right)|_{x_{2}=\cdots=x_{\ell}=x}.
\end{eqnarray}
Noticing that
$$P_{\ell+1}(x_{1},\dots,x_{\ell},x)=P_{\ell}(x_{2},\dots,x_{\ell},x)q(x_{1}/x)\prod_{j=2}^{\ell}q(x_{1}/x_{j}),$$
using (\ref{e-hypothesis}) and (\ref{eqab-comm}) we get
\begin{eqnarray}
q(x_{1}/x)^{\ell}Y_{W}(a,x_{1})Y_{W}(b,x)=q(x_{1}/x)^{\ell}Y_{W}(b,x)Y_{W}(a,x_{1}).
\end{eqnarray}
By Lemma \ref{lmodule-ab} we have
\begin{eqnarray*}
Y_{W}(a_{-1}b,x)=\left(q(x_{1}/x)^{\ell}Y_{W}(a,x_{1})Y_{W}(b,x)\right)|_{x_{1}=x}.
\end{eqnarray*}
Then we obtain
\begin{eqnarray*}
&&\left(P_{\ell+1}(x_{1},\dots,x_{\ell},x)
Y_{W}(a,x_{1})\cdots Y_{W}(a,x_{\ell})Y_{W}(a,x)\right)|_{x_{1}=\cdots=x_{\ell}=x}\\
&=&\left(q(x_{1}/x)^{\ell}Y_{W}(a,x_{1})P_{\ell}(x_{2},\dots,x_{\ell},x)Y_{W}(a,x_{2})\cdots Y_{W}(a,x_{\ell})Y_{W}(a,x)\right)|_{x_{1}=\cdots=x_{\ell}=x}\\
&=&\left(q(x_{1}/x)^{\ell}Y_{W}(a,x_{1})Y_{W}(b,x)\right)|_{x_{1}=x}\\
&=&Y_{W}(a_{-1}b,x)\\
&=&Y_{W}\left((a_{-1})^{\ell+1}{\bf 1},x\right),
\end{eqnarray*}
as desired.
\end{proof}

As an immediate consequence, we have (cf. \cite{dl}, \cite{ll},
\cite{dm}):

\bc{cnilpotent} Let $V$ be a vertex $\Gamma$-algebra, let $a\in V$
such that $a_{n}a=0$ for $n\ge 0$, and let $\ell$ be a positive
integer. Let $(W,Y_{W})$ be a $\Gamma$-equivariant quasi $V$-module
and let $q(z)$ be a polynomial as in Proposition \ref{pnilpotent}.
If $(a_{-1})^{\ell+1}{\bf 1}=0$ in $V$, then
\begin{eqnarray*}
&&\left(\prod_{1\le i<j\le \ell+1}q(x_{i}/x_{j})\right)
Y_{W}(a,x_{1})\cdots Y_{W}(a,x_{\ell})Y_{W}(a,x_{\ell+1})\nonumber\\
&\in& \Hom (W,W((x_{1},x_{2},\dots,x_{\ell+1})))
\end{eqnarray*}
and
\begin{eqnarray}\label{ePnilpotent}
\left(\left(\prod_{1\le i<j\le \ell+1}q(x_{i}/x_{j})\right)
Y_{W}(a,x_{1})\cdots Y_{W}(a,x_{\ell})Y_{W}(a,x_{\ell+1})\right)|_{x_{1}=\cdots =x_{\ell}=x_{\ell+1}}=0.
\end{eqnarray}
On the other hand, if (\ref{ePnilpotent}) holds and if
$(W,Y_{W})$ is faithful, then  $(a_{-1})^{\ell+1}{\bf 1}=0$.
\ec

Now, we consider the vertex $\Gamma_{q}$-algebra
$L_{\widehat{\A}}(\ell,0)$, recalling that $\Gamma_{q}=\{ q^{n}\ |\ n\in \Z\}$.

\bp{pmain-end} For any restricted $\hat{A}_{\hbar}$-module $W$ and
for any $\alpha\in \Z$,
$$\left(\prod_{1\le i<j\le k+1}(x_{i}/x_{j}-q^{2\alpha})(x_{i}/x_{j}-q^{-2\alpha})\right)
A_{\alpha}(x_{1})\cdots A_{\alpha}(x_{k+1})$$ lies in $\Hom
(W,W((x_{1},x_{2},\dots,x_{k+1})))$ for any positive integer $k$.
Furthermore, for a positive integer $\ell$, $\Gamma_{q}$-equivariant
quasi $L_{\widehat{\A}}(\ell,0)$-modules are exactly those
restricted $\hat{A}_{\hbar}$-modules $W$ of level $\ell$ such that
for every nonzero $\alpha\in \Z$,
$$\left(\prod_{1\le i<j\le \ell+1}(x_{i}/x_{j}-q^{2\alpha})(x_{i}/x_{j}-q^{-2\alpha})\right)
A_{\alpha}(x_{1})\cdots A_{\alpha}(x_{\ell+1})$$  vanishes at
$x_{i}/x_{j}=1$ for $1\le i<j\le \ell+1$. \ep

\begin{proof}  {}From (\ref{exalpha-xbeta}) we have
\begin{eqnarray}
(x-q^{\alpha+\beta}z)(x-q^{-(\alpha+\beta)}z)[A_{\alpha}(x),A_{\beta}(z)]=0
\end{eqnarray}
for $\alpha,\beta\in \Z$ (regardless whether $\alpha+\beta=0$). In
particular, we have
$$(x/z-q^{2\alpha})(x/z-q^{-2\alpha})[A_{\alpha}(x),A_{\alpha}(z)]=0. $$
Then the first assertion follows. As for the second assertion, in
view of Proposition \ref{pquotient-L}, $\Gamma_{q}$-equivariant
quasi $L_{\widehat{\A}}(\ell,0)$-modules are exactly those
$\Gamma_{q}$-equivariant quasi $V_{\widehat{\A}}(\ell,0)$-modules
$W$ on which $Y_{W}(G_{\alpha,m}(-1)^{\ell+1}{\bf 1},x)=0$ for all
$\alpha,m\in \Z$ with $\alpha\ne 0$. On the other hand, by Theorem
\ref{ttypeA-quasimodule}, $\Gamma_{q}$-equivariant quasi
$V_{\widehat{\A}}(\ell,0)$-modules are restricted
$\widehat{\A}_{\hbar}$-modules of level $\ell$. For $\alpha,m\in \Z$
with $\alpha\ne 0$, we have
$$(G_{\alpha,m})_{j}G_{\alpha,m}=0\ \ \mbox{ in }V_{\widehat{\A}}(\ell,0)
\ \ \mbox{ for }j\ge 0.$$ Then the second assertion follows
immediately from Corollary \ref{cnilpotent}.
\end{proof}

Define an anti-involution $\omega$ of $\widehat{A}_{\hbar}$ by
\begin{eqnarray}
\omega({\bf c})={\bf c},\ \ \ \ \omega(A_{\alpha,m})=A_{\alpha,-m}\
\ \ \mbox{ for }\alpha,m\in \Z.
\end{eqnarray}
From \cite{kr1}, a quasifinite module $L(\ell,\lambda)$ is unitary
if and only if there are  finitely many positive integers  $n_{i}$
and real numbers $\mu_{i}$  such that $\ell=\sum_{i}n_{i}$ and
 $$\lambda (A_{n,0})=
\frac{q^{n}}{q^{n}-q^{-n}}\sum_{i}n_{i}e^{n\mu_{i}}$$ for every
nonzero integer $n$. In particular, unitarity implies that $\ell$ is
a nonnegative integer. If $L(\ell, \lambda)$ is a unitary
quasifinite module of level $1$, then
$$\lambda (A_{n,0})=e^{n\mu}\cdot \frac{q^{n}}{q^{n}-q^{-n}}$$
for $n\in \Z\backslash \{0\}$, where $\mu$ is a real number.

\br{rD-Aconnection} {\em Recall from \cite{kr1} the Lie algebra
$\widehat{\D_{q}}$, which has a basis consisting of $C$ (a central
element) and $T_{m,n}$ $(m,n\in \Z)$, with commutator relations
$$[T_{m,n},T_{m',n'}]=2\sinh (h(m'n-mn')) T_{m+m',n+n'}+m\delta_{m+m',0}\delta_{n+n',0}C$$
for $m,n,m',n'\in \Z$, where $q=e^{2h}$. Lie algebras
$\widehat{\D_{q}}$ and $\widehat{A}_{\hbar}$ are isomorphic, where an
isomorphism from $\widehat{\D_{q}}$ to $\widehat{A}_{\hbar}$ is given by
$$C\mapsto {\bf c},\ \ \ \ T_{m,n}\mapsto -A_{n,m}
\ \ \mbox{ for }m,n\in \Z,$$ and by identifying  $e^{h}$ with
$e^{i\hbar}$.  For a weight $\lambda\in H^{*}$, set
$$\Delta_{\lambda}(x)=\sum_{n\ne 0}\Delta_{n}(\lambda)x^{-n},$$
where  for $n\ne 0$,
$$\Delta_{n}(\lambda)=\lambda(T_{0,n})=-\lambda(A_{n,0}).$$
 It was proved in \cite{kr1} that an irreducible highest weight $\widehat{\D_{q}}$-module
 $L(c,\lambda)$ is quasifinite if
and only if there exists a nonzero polynomial $b(x)$ such that
$$b(x)(\Delta_{\lambda}(x)-\Delta_{\lambda}(q^{-1}x)+c)=0.$$
Furthermore, it was proved that $L(c,\lambda)$ is unitary if and
only if there exist finitely many positive integers $n_{i}$ and real
numbers $a_{i}$ such that $c=\sum_{i}n_{i}$ and
$$\lambda(T_{0,n})=\sum_{i}\frac{n_{i}q^{na_{i}}}{1-q^{n}}
\ \ \ \mbox{ for all }n\in \Z\backslash \{0\}.$$
Correspondingly,
for $\widehat{A}_{\hbar}$ we have
$$\lambda(A_{n,0})=-\sum_{i}\frac{n_{i}q^{2na_{i}}}{1-q^{2n}}
=\sum_{i}n_{i}q^{2n(a_{i}-1)}\frac{q^{n}}{q^{n}-q^{-n}},$$ (noticing
that $q=e^{2h}$ for $\widehat{\D_{q}}$ whereas $q=e^{i\hbar}$ for
$\widehat{A}_{\hbar}$).} \er

As the main result of this section we have:

\bt{tmain-classification} Let $\ell$ be a positive integer. Then
every unitary quasifinite highest weight irreducible
$\widehat{A}_{\hbar}$-module of level $\ell$ is an irreducible
$\Gamma_{q}$-equivariant quasi $L_{\widehat{\A}}(\ell,0)$-module.
\et

\begin{proof} We here shall apply Proposition \ref{pmain-end}.
First we consider the case $\ell=1$. Set
$$M(1)=\C[x_{1},x_{2},x_{3},\dots],$$
the algebra of polynomials in variables $x_{1},x_{2},\dots$. Let
$\mu\in \R$. For $\alpha\in \Z\backslash \{0\}$, set
\begin{eqnarray*}
\widehat{X}_{\alpha}(z)=a_{\alpha}\exp\left(\sum_{m\ge
1}z^{m}(q^{m\alpha}-q^{-m\alpha})x_{m}\right)\exp\left(\sum_{m\ge
1}\frac{z^{-m}}{m}(q^{m\alpha}-q^{-m\alpha})\frac{\partial}{\partial
x_{m}}\right),
\end{eqnarray*}
 where $q=e^{i\hbar}$ as before and
$a_{\alpha}=e^{i\mu\alpha}\frac{q^{\alpha}}{q^{\alpha}-q^{-\alpha}}\in
\C$. It was proved in \cite{gkl} that $M(1)$ becomes an irreducible
highest weight $\widehat{A}_{\hbar}$-module of level $1$ with
\begin{eqnarray}
&&A_{0,m}=\frac{\partial}{\partial x_{m}},\ \ \ A_{0,-m}=mx_{m}\ \ \
\mbox{ for }m\ge 1,\nonumber\\
 &&A_{\alpha}(z)=\widehat{X}_{\alpha}(z)\ \
\mbox{ for }\alpha\in \Z\backslash \{0\},
\end{eqnarray}
where the highest weight $\lambda$ is determined by $\lambda_{0}=0$
and $\lambda_{\alpha}=a_{\alpha}$ for $\alpha\in \Z\backslash
\{0\}$. Denote this $\widehat{A}_{\hbar}$-module by $M(1)^{[\mu]}$.
 From \cite{gkl}, for
$\alpha,\beta\in \Z\backslash \{0\}$ we have
\begin{eqnarray}
\widehat{X}_{\alpha}(z)\widehat{X}_{\beta}(w)=\frac{(z-q^{\alpha-\beta}w)(z-q^{\beta-\alpha}w)}
{(z-q^{\alpha+\beta}w)(z-q^{-(\alpha+\beta)}w)}\NO
\widehat{X}_{\alpha}(z)\widehat{X}_{\beta}(w)\NO.
\end{eqnarray}
That is,
\begin{eqnarray*}
&&\left(\frac{z}{w}-q^{\alpha+\beta}\right)\left(\frac{z}{w}-q^{-(\alpha+\beta)}\right)
\widehat{X}_{\alpha}(z)\widehat{X}_{\beta}(w)\nonumber\\
&=&\left(\frac{z}{w}-q^{\alpha-\beta}\right)
\left(\frac{z}{w}-q^{\beta-\alpha}\right)\NO
\widehat{X}_{\alpha}(z)\widehat{X}_{\beta}(w)\NO.
\end{eqnarray*}
{}From this we have
\begin{eqnarray}
\left(1-\frac{q^{2\alpha}w}{z}\right)\left(1-\frac{q^{-2\alpha}w}{z}\right)
\widehat{X}_{\alpha}(z)\widehat{X}_{\alpha}(w)\in \Hom (W,W((z,w)))
\end{eqnarray}
and
\begin{eqnarray}
\left(\left(1-\frac{q^{2\alpha}w}{z}\right)\left(1-\frac{q^{-2\alpha}w}{z}\right)
\widehat{X}_{\alpha}(z)\widehat{X}_{\alpha}(w)\right)|_{z=w}=0.
\end{eqnarray}
Furthermore, for $m\in \Z$, we have
\begin{eqnarray}
\left(1-q^{2\alpha}\frac{w}{z}\right)\left(1-q^{-2\alpha}\frac{w}{z}\right)
\widehat{X}_{\alpha}(q^{m}z)\widehat{X}_{\alpha}(q^{m}w)\in \Hom (W,W((z,w)))
\end{eqnarray}
and it vanishes at $z=w$. It follows from Proposition
\ref{pmain-end} that $M(1)^{[\mu]}$ is an irreducible
$\Gamma_{q}$-equivariant quasi $L_{\widehat{\A}}(1,0)$-module. On
the other hand, it was proved in \cite{kr1} that every unitary
quasifinite highest weight irreducible $\widehat{A}_{\hbar}$-module of
level $1$ is of this form. This proves the theorem for $\ell=1$.

As for the high level case (with $\ell>1$), we apply the common
trick---to use tensor products. Let $L(\ell,\lambda)$ be a unitary
quasifinite irreducible highest weight $\widehat{A}_{\hbar}$-module of
level $\ell$. Then from Remark \ref{rD-Aconnection} there exist
finitely many positive integers $n_{j}$ and real numbers $\mu_{j}$
such that $\ell=\sum_{j}n_{j}$ and
$$\lambda(A_{n,0})=\frac{q^{n}}{q^{n}-q^{-n}}\sum_{j}n_{j}e^{in\mu_{j}}$$
for all $n\in \Z\backslash \{0\}$. Consider the tensor product
$\widehat{A}_{\hbar}$-module
$$M=\otimes_{j} \left(M(1)^{[\mu_{j}]}\right)^{\otimes n_{j}},$$
which is a restricted $\widehat{A}_{\hbar}$-module of level $\ell$. For
every nonzero $\alpha\in \Z$,
$$\left(\prod_{1\le i<j\le \ell+1}(x_{i}/x_{j}-q^{2\alpha})(x_{i}/x_{j}-q^{-2\alpha})\right)
A_{\alpha}(x_{1})\cdots A_{\alpha}(x_{\ell+1})$$ vanishes  on $M$ at
$x_{i}/x_{j}=1$ for $1\le i<j\le \ell+1$. Then it follows from
Proposition \ref{pmain-end} that $M$ is a $\Gamma_{q}$-equivariant
quasi $L_{\widehat{\A}}(\ell,0)$-module. Set
$v=\otimes_{j}1^{\otimes n_{j}}\in M$. We see that $v$ is a highest
weight vector of weight $\lambda$. It follows from unitarity that
the submodule $U(\widehat{A}_{\hbar})v$ generated by $v$ is an
irreducible module. Consequently, $U(\widehat{A}_{\hbar})v$ is
isomorphic to $L(\ell,\lambda)$. Thus $L(\ell,\lambda)$ is a
$\Gamma_{q}$-equivariant quasi $L_{\widehat{\A}}(\ell,0)$-module.
\end{proof}

We end up this paper with the following:

\begin{conj}
Every $\N$-graded irreducible
$\Gamma_{q}$-equivariant quasi $L_{\widehat{\A}}(\ell,0)$-module is
a unitary quasifinite highest weight irreducible
$\widehat{A}_{\hbar}$-module of level $\ell$.
\end{conj}


\begin{thebibliography}{AAGBP}



\bibitem[DL]{dl}
C. Dong and J. Lepowsky, {\em Generalized Vertex Algebras and
Relative Vertex Operators}, Progress in Math., Vol. {\bf 112},
Birkh\"auser, Boston, 1993.

\bibitem[DM]{dm}
J. Ding and Miwa, Zeros and poles of quantum current operators and
the condition of quantum integrability, arXiv:q-alg/9608001.

\bibitem[DF]{df}
J. Ding and B. Feigin, Quantum current operators (III): Commutative
quantum current operators, semi-infinite construction and functional
models, arXiv: q-alg/9612009.

\bibitem[DLM]{dlm}
C. Dong, H.-S. Li and G. Mason, Regularity of rational vertex
operator algebras, {\em Adv. Math.} {\bf 132} (1997) 148-166

\bibitem[FFZ]{ffz}
D. Fairlie, P. Fletcher, C. Zachos, Trigonometric structure
constants for new infinite-dimensional algebras, {\em Phys. Lett.}
{\bf B 218} (1989) 203-206.

\bibitem[F]{flor}
E. G. Floratos, Spin wedge and vertex operator representations of
trigonometric algebras and their central extensions, {\em Phys.
Lett.} {\bf B 232} (1989) 467-474.

\bibitem[FLM]{flm}
I. B. Frenkel, J. Lepowsky and A. Meurman, {\em Vertex Operator
Algebras and the Monster,} Pure and Applied Math., Vol. 134,
Academic Press, Boston, 1988.

\bibitem[FZ]{fz}
I. B. Frenkel and Y.-C. Zhu,  Vertex operator algebras associated to
representations of affine and Virasoro algebras, {\em Duke Math. J.}
{\bf 66} (1992) 123-168.

\bibitem[G-KK]{gkk}
M. Golenishcheva-Kutuzova and V. Kac, $\Gamma$-conformal algebras,
{\em J. Math. Phys.} {\bf 39} (1998) 2290-2305.

\bibitem[G-KL1]{gkl1}
M. Golenishcheva-Kutuzova and D. Lebedev, $\Z$-graded trigonometric
Lie subalgebras in $\hat{A}_{\infty}$, $\hat{B}_{\infty}$,
$\hat{C}_{\infty}$, and $\hat{D}_{\infty}$ and their vertex operator
representations, {\em Funktsional'nyi Analiz i Ego Prilozhenniya}
{\bf 27} (1993) 12-24.

\bibitem[G-KL2]{gkl}
M. Golenishcheva-Kutuzova and D. Lebedev, Vertex operator
representation of some quantum tori Lie algebras, {\em Commun. Math.
Phys.} {\bf 148} (1992) 403-416.

\bibitem [H]{hoppe}
J. Hoppe, ${\rm Diff}_{A}T^{2}$ and the curvature of some
infinite-dimensional manifolds, {\em Phys. Lett.} {\bf B 215} (1988)
706-710.

\bibitem [K]{kac1}
V. G. Kac, {\it Infinite-dimensional Lie Algebras}, 3rd ed.,
Cambridge Univ. Press, Cambridge, 1990.

\bibitem [KR]{kr1}
V. G. Kac and A. Radul, Quasifinite highest weight modules over the
Lie algebra of differential operators on the circle,  {\em Commun.
Math. Phys.} {\bf 157} (1993) 429-457.

\bibitem[LL]{ll}
J. Lepowsky and H.-S. Li, {\em Introduction to Vertex Operator
Algebras and Their Representations}, Progress in Math. {\bf 227},
Birkh\"auser, Boston, 2004.

\bibitem[LP]{lp}
J. Lepowsky and M. Primc, {\em Structure of the Standard Modules for
the Affine Lie Algebra $A_{1}^{(1)}$,} Contemporary Math. {\bf 46},
Amer. Math. Soc., Providence, 1985.

\bibitem[Li1]{li-local}
H.-S. Li,  Local systems of vertex operators, vertex superalgebras
and modules,  {\em J. Pure Appl. Algebra} {\bf 109} (1996) 143-195.

\bibitem[Li2]{li-twisted}
H.-S. Li,  Local systems of twisted vertex operators, vertex superalgebras
and twisted modules, Contemporary Math. {\bf 193},
Amer. Math. Soc., Providence, 1996, 203-236.

\bibitem[Li3]{li-gamma}
H.-S. Li, A new construction of vertex algebras and quasi modules
for vertex algebras, {\em Adv. Math.} {\bf 202} (2006) 232-286.

\bibitem[Li4]{li-tlie}
H.-S. Li, On certain generalizations of twisted affine Lie algebras
and quasimodules for $\Gamma$-vertex algebras, {\em J. Pure Appl.
Algebra} {\bf 209} (2007) 853-871.

\bibitem[Li5]{li-twisted-quasi}
H.-S. Li, Twisted modules and quasi-modules for vertex operator
algebras, Contemporary Math. {\bf 422}, Amer. Math. Soc.,
Providence, 2007, 389-400.

\bibitem[LTW]{ltw-twisted}
H.-S. Li, S.-B. Tan and Q. Wang, Twisted modules for quantum vertex algebras,
{\em J. Pure Appl. Algebra} {\bf 214} (2010) 201-220.

\bibitem[MP1]{mp1}
A. Meurman and M. Primc, Vertex operator algebras and representations of affine Lie algebras,
{\em Acta Applicandae Math.} {\bf 44} (1996), 207-215.

\bibitem[MP2]{mp}
A. Meurman and M. Primc, {\em Annihilating Fields of Standard Modules of $\widetilde{\sl(2,\C)}$ and
Combinatorial Identities}, preprint 1994; Memoirs Amer. Math. Soc. {\bf 652}, 1999.

\end{thebibliography}
\end{document}